\documentclass[11pt]{article}
\usepackage{}
\usepackage{amssymb}
\usepackage{amsfonts}
\usepackage{mathrsfs}
\usepackage{graphicx}
\usepackage{amsmath}
\usepackage{color}
\usepackage{amsthm}

\usepackage{pifont,bm}
\usepackage{tikz}
\usepackage{enumerate}
\theoremstyle{definition}

\makeatletter
\def\@biblabel#1{[#1]}
\makeatother

\makeatletter \@addtoreset{equation}{section}

\begin{document}
%\begin{CJK*}{GBK}{song}

\begin{titlepage}
\title{\bf{Riemann-Hilbert problem of the three-component coupled Sasa-Satsuma equation and its multi-soliton solutions
\footnote{This work is supported by the National Key Research and Development
Program of China under Grant No. 2017YFB0202901 and the National Natural Science Foundation of China under Grant No.11871180.\protect\\
\hspace*{3ex}$^{*}$Corresponding authors.\protect\\
\hspace*{3ex} E-mail addresses: xiubinwang@163.com (X.B. Wang), bohan@hit.edu.cn (B. Han)}
}}
\author{Xiu-Bin Wang$^{*}$, Bo Han \\
%%%%%%%%%%%%%%%%%%%%%%%%%%%%%%%%%%%%%%%%%%%%%%%%%%%%%%%%%%%%%%%%%%%%%%%%%%%%%%%%%%%%%%%%%
%%%%%              以下两行为作者单位
%%%%%%%%%%%%%%%%%%%%%%%%%%%%%%%%%%%%%%%%%%%%%%%%%%%%%%%%%%%%%%%%%%%%%%%%%%%%%%%%%%%%%%%%%
\small \emph{Department of Mathematics, Harbin Institute of Technology, Harbin 150001, People's Republic of China}\\
%\small \emph{$^{b}$School of Mathematics and Institute of Mathematical Physics, China University of Mining and Technology,} \\
%\small \emph{Xuzhou 221116, People's Republic of China}\\
\date{}}
\thispagestyle{empty}
\end{titlepage}
\maketitle

\vspace{-0.5cm}
\begin{center}
\rule{15cm}{1pt}\vspace{0.3cm}

\parbox{15cm}{\small
{\bf Abstract}\\
\hspace{0.5cm} In this work,
the inverse scattering transform of the three-component coupled Sasa-Satsuma equation is investigated via the Riemann-Hilbert method.
Firstly we consider a Lax pair associated with a $7\times 7$ matrix spectral problem for the equation.
Then we present the spectral analysis of the Lax pair, from which  a kind of Riemann-Hilbert problem is formulated.
Moreover, $N$-soliton solutions to the equation are constructed through a particular Riemann-Hilbert problem with vanishing scattering coefficients. Finally, the dynamics of the soliton solutions are discussed with some graphics.}

\vspace{0.5cm}
\parbox{15cm}{\small{

\vspace{0.3cm} \emph{Key words:} The three-component coupled Sasa-Satsuma equation; Riemann-Hilbert problem; Multi-soliton solutions.

 } }
\end{center}
\vspace{0.3cm} \rule{15cm}{1pt} \vspace{0.2cm}

\section{Introduction}

It is well known that much physical phenomena can be proved or obtained by the soliton solutions of nonlinear evolutions equations (NLEEs).
As a consequence, the investigation of soliton solutions for NLEEs equations has become more and more attractive.
A variety of approaches have been proposed over these years for seeking soliton solutions \cite{vbm-1991}-\cite{RH-2004}.
Among those methods,
the Riemann-Hilbert (RH) approach is one of the most effective methods to seek soliton solutions for integrable systems.
The main process of this method is to seek
a corresponding Riemann-Hilbert problem (RHP) on the spectral analysis of integrable systems.
In 1980s, Beals and his co-workers \cite{mmmm-1} investigated direct and inverse scattering for Ablowitz-Kaup-Newell-Segur systems on the
line with integrable matrix-valued potentials and, derived partial characterization of the scattering data.
These works are very familiar with the various boundary conditions for the Jost function,
under the condition of the location of the spectral parameter.
More importantly, it is very helpful for us to construct the RHP of integrable systems.
Recently, the RH approach is also extended to consider initial boundary
value problems and asymptotics of integrable equations \cite{sf-1}-\cite{as-1}.

The standard nonlinear Schr\"{o}dinger (NLS) equation has been paid much attention due to its widespread applications in optics, Bose-Einstein condensates, hydrodynamics, plasma physics, molecular biology and even finance.
To well describe other important types of nonlinear physical phenomena in a similar way, it is necessary to go beyond the standard
NLS description. One prime research is to add higher-order terms and/or dissipative terms to the
NLS equation to accurately model extreme wave events in some nonlinear wave systems such as micro-structured optical fibres and fibre lasers.
Another important development focuses on
the investigation of coupled-wave systems, as many physical systems comprise interacting wave
components of distinct modes, frequencies or polarizations. Recently, the coupled NLS
equations have become a topic of intense research, since the components are usually more than
one practically for many physical phenomena.
%To the best knowledge of the authors, the RH method has been used to consider some crucial integrable equations, such as
%the vector NLS equation, the SS equation,
%the general coupled NLS equation etc
%Because of the RHP and soliton solutions for multi-component NLS equations are more interesting and complex than scalar NLS equation,
Therefore, the chief aim of this work is to construct the RHP and soliton solutions
for the following three-coupled higher-order NLS equation, whose form reads \cite{kkpa-1998}
\begin{equation}\label{mc-CNLS}
\left\{ \begin{aligned}
&iq_{1T}+\frac{1}{2}q_{1XX}+q_{1}\sum_{N=1}^{3}|q_{N}|^{2}+i\left[q_{1XXX}+6q_{1X}\sum_{N=1}^{3}|q_{N}|^{2}\right.\\
&\left.+3q_{1}\left(\sum_{N=1}^{3}|q_{N}|^{2}\right)_{X}\right]=0,\\
&iq_{2T}+\frac{1}{2}q_{2XX}+q_{2}\sum_{N=1}^{3}|q_{N}|^{2}+i\left[q_{2XXX}+6q_{2X}\sum_{N=1}^{3}|q_{N}|^{2}\right.\\
&\left.+3q_{2}\left(\sum_{N=1}^{3}|q_{N}|^{2}\right)_{X}\right]=0,\\
&iq_{3T}+\frac{1}{2}q_{3XX}+q_{3}\sum_{N=1}^{3}|q_{N}|^{2}+i\left[q_{3XXX}+6q_{3X}\sum_{N=1}^{3}|q_{N}|^{2}\right.\\
&\left.+3q_{3}\left(\sum_{N=1}^{3}|q_{N}|^{2}\right)_{X}\right]=0,
                 \end{aligned} \right.
\end{equation}
where $q_{1}=q_{1}(X,T)$, $q_{2}=q_{2}(X,T)$ and $q_{3}=q_{3}(X,T)$  are three complex functions of (X,T).
As mentioned in Ref.\cite{kkpa-1998}, in order to investigate Eq.\eqref{mc-CNLS} conveniently,
it is very necessary to rewrite \eqref{mc-CNLS} in the three-component coupled Sasa-Satsuma equation
\begin{equation}\label{SSE}
\left\{ \begin{aligned}
u_{1,t}+&u_{1,xxx}+6\left(|u_{1}|^{2}+|u_{2}|^{2}+|u_{3}|^{2}\right)u_{1,x}+3u_{1}\left(|u_{1}|^{2}+|u_{2}|^{2}+|u_{3}|^{2}\right)_{x}=0,\\
u_{2,t}+&u_{2,xxx}+6\left(|u_{1}|^{2}+|u_{2}|^{2}+|u_{3}|^{2}\right)u_{2,x}+3u_{2}\left(|u_{1}|^{2}+|u_{2}|^{2}+|u_{3}|^{2}\right)_{x}=0,\\
u_{3,t}+&u_{3,xxx}+6\left(|u_{1}|^{2}+|u_{2}|^{2}+|u_{3}|^{2}\right)u_{3,x}+3u_{3}\left(|u_{1}|^{2}+|u_{2}|^{2}+|u_{3}|^{2}\right)_{x}=0,
                 \end{aligned} \right.
\end{equation}
by making use of the following three transformations (i.e., gauge, Galilean and scale transformations)
\begin{equation}\label{FT-1}
\left\{ \begin{aligned}
&u_{1}(x,t)=q_{1}(X,T)\exp\left\{-\frac{i}{6}\left(X-\frac{T}{18}\right)\right\},\\
&u_{2}(x,t)=q_{2}(X,T)\exp\left\{-\frac{i}{6}\left(X-\frac{T}{18}\right)\right\},\\
&u_{3}(x,t)=q_{3}(X,T)\exp\left\{-\frac{i}{6}\left(X-\frac{T}{18}\right)\right\},\\
&x=X-\frac{T}{12},~~t=T.
            \end{aligned} \right.
\end{equation}

%Much research has been done for the special cases of the Eq. \eqref{SSE},
%but the RHP of Eq.\eqref{SSE} has not been studied so far.
As we all know that the RH approach is a effective way to
construct soliton solutions \cite{wx-1}-\cite{yjk-2003}. Nonetheless, since the Eq. \eqref{SSE} involves
a $7\times 7$ matrix spectral problem, the RHP for the Eq. \eqref{SSE} is rather hard to deal
with. The research in this direction, to our best knowledge, has
not been considered so far. The chief purpose of the present article
is to discuss the RHP and soliton solutions of the Eq. \eqref{SSE} by utilizing
the RH approach.

The layout of this paper is arranged as follows.
In section II, we present the spectral analysis of the Eq. \eqref{SSE},
from which a kind of RHP is established.
In section III, by considering a specific RHP, we obtain multi-soliton solutions of the Eq. \eqref{SSE}.
Besides, some figures are presented to understand the dynamics of these solutions.
The last section summarizes the results of this work.

\section{Riemann-Hilbert problem}
The Eq. \eqref{SSE} admits the following Lax pair \cite{kkpa-1998}
\begin{equation}\label{Lax-1}
\left\{ \begin{aligned}
&\Phi_{x}=U\Psi,~~U=i\lambda\sigma_{3}+Q,\\
&\Phi_{t}=V\Psi,~~V=4i\lambda^3\sigma_{3}+4\lambda^{2}Q+2i\lambda\left(Q^2+Q_{x}\right)\sigma_{3}
+Q_{x}Q-QQ_{x}-Q_{xx}+2Q^3,
            \end{aligned} \right.
\end{equation}
with
\begin{equation}\label{Lax-2}
Q=\left(
    \begin{array}{ccccccc}
      0 & 0 & 0 & 0 & 0 & 0 & u_{1} \\
      0 & 0 & 0 & 0 & 0 & 0 & \bar{u}_{1} \\
      0 & 0 & 0 & 0 & 0 & 0 & u_{2} \\
      0 & 0 & 0 & 0 & 0 & 0 & \bar{u}_{2} \\
      0 & 0 & 0 & 0 & 0 & 0 & u_{3} \\
      0 & 0 & 0 & 0 & 0 & 0 & \bar{u}_{3} \\
      -\bar{u}_{1} & -u_{1} & -\bar{u}_{2} & -u_{2} & -\bar{u}_{3} & -u_{3} & 0 \\
    \end{array}
  \right),
\end{equation}
where $\sigma_{3}=\mbox{diag}(1,1,1,1,1,1,-1)$, and $\Phi=\Phi(x,t,\lambda)$ is a column vector function of the spectral parameter $\lambda$.
Additionally,
the compatibility condition of the Lax pair \eqref{Lax-1} yields the Eq. \eqref{SSE}.

For convenience, we consider a new matrix spectral function $\Psi=\Psi(x,t,\lambda)$ given by
\begin{equation}\label{RHP-1}
\Phi(x,t,\lambda)=\Psi(x,t,\lambda) e^{i\left(\lambda x+4i\lambda^{3}t\right)\sigma_{3}}.
\end{equation}
Then the spectral problem \eqref{Lax-1} can be rewritten as
\begin{equation}\label{RHP-2}
\left\{ \begin{aligned}
&\Psi_{x}+i\lambda[\Psi,\sigma_{3}]=Q\Psi,\\
&\Psi_{t}+4i\lambda^3[\Psi,\sigma_{3}]=\widetilde{Q}\Psi,
             \end{aligned} \right.
\end{equation}
with
\begin{equation}\label{RHP-3}
\widetilde{Q}=4\lambda^{2}Q+2i\lambda\left(Q^2+Q_{x}\right)\sigma_{3}+Q_{x}Q-QQ_{x}-Q_{xx}+2Q^3.
\end{equation}

Throughout this article,
$\mathbb{C}^{+}=\{z\in\mathbb{C}|\Im(z)>0\}$,  $\mathbb{C}^{-}=\{z\in\mathbb{C}|\Im(z)<0\}$
represent the upper half-plane and the lower half-plane, respectively.
Let us next consider the two matrix functions
\begin{align}\label{RHP-4}
\Psi_{\pm}=
\left[(\Psi_{\pm})_{1},(\Psi_{\pm})_{2},(\Psi_{\pm})_{3},(\Psi_{\pm})_{4},(\Psi_{\pm})_{5},(\Psi_{\pm})_{6},(\Psi_{\pm})_{7}\right],
\end{align}
under the following asymptotic conditions
\begin{equation}\label{RHP-5}
\Psi_{+}\rightarrow\mathbb{I},~~~~x\rightarrow +\infty~~~~\mbox{and}~~~~\Psi_{-}\rightarrow\mathbb{I},~~~~x\rightarrow -\infty.
\end{equation}
Here each $(\Psi_{\pm})_{m}$ represents the mth column of $\Psi_{\pm}$, respectively.
The matrix $\mathbb{I}$ is the $7\times 7$ unit matrix, and
the subscripts of $\Psi$ refer to which end of the x-axis the boundary conditions are required.
Utilizing the above boundary condition \eqref{RHP-5},
the two matrix solutions $\Psi_{\pm}$ can be given by the following two Volterra integral equations
\begin{equation}\label{RHP-6}
\left\{ \begin{aligned}
&\Psi_{+}(x,\lambda)
=\mathbb{I}+\int_{+\infty}^{x}e^{i\lambda\sigma_{3}(x-y)}Q(y)\Psi_{+}(y,\lambda)e^{i\lambda\sigma_{3}(y-x)}dy,\\
&\Psi_{-}(x,\lambda)
=\mathbb{I}+\int_{-\infty}^{x}e^{i\lambda\sigma_{3}(x-y)}Q(y)\Psi_{-}(y,\lambda)e^{i\lambda\sigma_{3}(y-x)}dy,
                          \end{aligned} \right.
\end{equation}
respectively.
It is not hard to check that
\begin{equation}
(\Psi_{+})_{1},(\Psi_{+})_{2},(\Psi_{+})_{3},(\Psi_{+})_{4},(\Psi_{+})_{5},(\Psi_{+})_{6},(\Psi_{-})_{7},
\end{equation}
can be analytically extendible to $\mathbb{C}^{+}$. Besides,
\begin{equation}
(\Psi_{-})_{1},(\Psi_{-})_{2},(\Psi_{-})_{3},(\Psi_{-})_{4},(\Psi_{-})_{5},(\Psi_{-})_{6},(\Psi_{+})_{7},
\end{equation}
can be analytically extendible to $\mathbb{C}^{-}$.

In the following, we study the analytic properties of $\Psi_{\pm}$. Actually, noticing $\mbox{tr}(Q)=\mbox{tr}(\widetilde{Q})=0$,
we have $\det(\Psi_{\pm})=1,~~~~\lambda\in\mathbb{R}$.
%\begin{equation}\label{RHP-7}
%\det(\Psi_{\pm})=1,~~~~\lambda\in\mathbb{R}.
%\end{equation}
Besides, $\Psi_{\pm}$ can be linearly related by
\begin{equation}\label{RHP-8}
\Psi_{-}e^{i\lambda\sigma_{3} x}=\Psi_{+}e^{i\lambda\sigma_{3} x} \Omega(\lambda),~~\lambda\in\mathbb{R},
\end{equation}
where $\Omega(\lambda)=(\Omega_{ij})_{7\times7}, \lambda\in\mathbb{R}$.
%\begin{equation}\label{RHP-9}
%\Omega(\lambda)=\left(
%             \begin{array}{ccccccc}
%               \Omega_{11} & \Omega_{12} & \Omega_{13} & \Omega_{14} & \Omega_{15} & \Omega_{16} & \Omega_{17} \\
%               \Omega_{21} & \Omega_{22} & \Omega_{23} & \Omega_{24} & \Omega_{25} & \Omega_{26} & \Omega_{27} \\
%               \Omega_{31} & \Omega_{32} & \Omega_{33} & \Omega_{34} & \Omega_{35} & \Omega_{36} & \Omega_{37} \\
%               \Omega_{41} & \Omega_{42} & \Omega_{43} & \Omega_{44} & \Omega_{45} & \Omega_{46} & \Omega_{47} \\
%               \Omega_{51} & \Omega_{52} & \Omega_{53} & \Omega_{54} & \Omega_{55} & \Omega_{56} & \Omega_{57} \\
%               \Omega_{61} & \Omega_{62} & \Omega_{63} & \Omega_{64} & \Omega_{65} & \Omega_{66} & \Omega_{67} \\
%               \Omega_{71} & \Omega_{72} & \Omega_{73} & \Omega_{74} & \Omega_{75} & \Omega_{76} & \Omega_{77} \\
%             \end{array}
%           \right),~~\lambda\in\mathbb{R}.
%\end{equation}
Noticing that $\det(\Omega(\lambda))=1$ since $\det(\Psi_{\pm})=1$.
From the analytic property of $\Psi_{-}$, we know that $\Omega_{77}$ can be analytically extended to $\mathbb{C}^{+}$,
otherwise $\Omega_{kj} (1\leq k,j\leq 6)$ can be analytically extended to $\mathbb{C}^{-}$.
Generally, $\Omega_{k7}$, $\Omega_{7j}$ are not extended off the real $\lambda$-axis.

By considering the analytic properties of $\Psi_{\pm}$, we can obtain
\begin{align}\label{RHP-10}
P_{1}=
\left[(\Psi_{+})_{1},(\Psi_{+})_{2},(\Psi_{+})_{3},(\Psi_{+})_{4},(\Psi_{+})_{5},(\Psi_{+})_{6},(\Psi_{-})_{7}\right],
\end{align}
which is analytic in $\mu\in\mathbb{C}^{+}$. In addition, we have the following asymptotic behavior of $P_{1}$ (i.e.,$P_{1}(\lambda)\rightarrow\mathbb{I},~~~~\lambda\rightarrow\infty.$)
%\begin{equation}\label{RHP-11}
%P_{1}(\lambda)\rightarrow\mathbb{I},~~~~\lambda\rightarrow\infty.
%\end{equation}
In order to further derive RHP for the Eq. \eqref{SSE},
we must construct a analytic matrix $P_{2}$ in $\mathbb{C}^{-}$.
For this purpose, we first introduce the adjoint equation of the first expression in \eqref{Lax-1}
\begin{equation}\label{RHP-12}
K_{x}+i\lambda[K,\sigma_{3}]=QK.
\end{equation}
It is not hard to know that $\Psi_{\pm}^{-1}$ meet the above expression \eqref{RHP-12}.
Then let us introduce
%\begin{equation}\label{RHP-13}
%\left\{ \begin{aligned}
%  &\Psi_{+}^{-1}=\left[(\Psi_{+}^{-1})^{1},(\Psi_{+}^{-1})^{2},(\Psi_{+}^{-1})^{3},(\Psi_{+}^{-1})^{4},
%  (\Psi_{+}^{-1})^{5},(\Psi_{+}^{-1})^{6},(\Psi_{+}^{-1})^{7}\right]^{T},\\
%  &\Psi_{-}^{-1}=\left[(\Psi_{-}^{-1})^{1},(\Psi_{-}^{-1})^{2},(\Psi_{-}^{-1})^{3},(\Psi_{-}^{-1})^{4},
%  (\Psi_{-}^{-1})^{5},(\Psi_{-}^{-1})^{6},(\Psi_{-}^{-1})^{7}\right]^{T}.
%     \end{aligned} \right.
%\end{equation}
\begin{equation}\label{RHP-13}
\Psi_{+}^{-1}=\left(
        \begin{array}{c}
          (\Psi_{+}^{-1})^{1} \\
          (\Psi_{+}^{-1})^{2} \\
          (\Psi_{+}^{-1})^{3} \\
          (\Psi_{+}^{-1})^{4} \\
          (\Psi_{+}^{-1})^{5} \\
          (\Psi_{+}^{-1})^{6} \\
          (\Psi_{+}^{-1})^{7} \\
        \end{array}
      \right),~~
\Psi_{-}^{-1}=\left(
        \begin{array}{c}
          (\Psi_{-}^{-1})^{1} \\
          (\Psi_{-}^{-1})^{2} \\
          (\Psi_{-}^{-1})^{3} \\
          (\Psi_{-}^{-1})^{4} \\
          (\Psi_{-}^{-1})^{5} \\
          (\Psi_{-}^{-1})^{6} \\
          (\Psi_{-}^{-1})^{7} \\
        \end{array}
      \right).
\end{equation}
Here each $(\Psi_{\pm}^{-1})^{m}$ represents the mth row of $\Psi_{\pm}^{-1}$, respectively.
Hence we can find that
\begin{equation}\label{RHP-14}
\left\{ \begin{aligned}
  &\left[(\Psi_{+}^{-1})^{1},(\Psi_{+}^{-1})^{2},(\Psi_{+}^{-1})^{3},(\Psi_{+}^{-1})^{4},
  (\Psi_{+}^{-1})^{5},(\Psi_{+}^{-1})^{6},(\Psi_{-}^{-1})^{7}\right],\\
  &\left[(\Psi_{-}^{-1})^{1},(\Psi_{-}^{-1})^{2},(\Psi_{-}^{-1})^{3},(\Psi_{-}^{-1})^{4},
  (\Psi_{-}^{-1})^{5},(\Psi_{-}^{-1})^{6},(\Psi_{+}^{-1})^{7}\right],
             \end{aligned} \right.
\end{equation}
can be analytically extended to $\mathbb{C}^{-}$ and $\mathbb{C}^{+}$, respectively.
As a result, we can give a matrix function
\begin{equation}\label{RHP-15}
\left(
        \begin{array}{c}
          (\Psi_{+}^{-1})^{1} \\
          (\Psi_{+}^{-1})^{2} \\
          (\Psi_{+}^{-1})^{3} \\
          (\Psi_{+}^{-1})^{4} \\
          (\Psi_{+}^{-1})^{5} \\
          (\Psi_{+}^{-1})^{6} \\
          (\Psi_{-}^{-1})^{7} \\
        \end{array}
      \right),
\end{equation}
which is analytic in $\lambda\in\mathbb{C}^{-}$, and the asymptotic behavior of $P_{2}$ yields $P_{2}(\lambda)\rightarrow\mathbb{I},~~\lambda\rightarrow\infty.$
%\begin{equation}\label{RHP-16}
%P_{2}(\lambda)\rightarrow\mathbb{I},~~\lambda\rightarrow\infty.
%\end{equation}
Moreover, we can easily find that $e^{-i\lambda\sigma_{3}x}\Psi_{\pm}^{-1}$
can be linearly related by a scattering matrix $\Theta(\lambda)=\Omega^{-1}(\lambda)$
\begin{equation}\label{RHP-17}
e^{-i\lambda\sigma_{3}x}\Psi_{-}^{-1}=\Theta(\lambda)e^{-i\lambda\sigma_{3}x}\Psi_{+}^{-1},~~\lambda\in\mathbb{R}.
\end{equation}
where %$\Omega(\lambda)=(\Theta_{ij})_{7\times7}, \lambda\in\mathbb{R}$.
\begin{equation}\label{RHP-18}
\Theta(\lambda)=\left(
             \begin{array}{ccccccc}
               \Theta_{11} & \Theta_{12} & \Theta_{13} & \Theta_{14} & \Theta_{15} & \Theta_{16} & \Theta_{17} \\
               \Theta_{21} & \Theta_{22} & \Theta_{23} & \Theta_{24} & \Theta_{25} & \Theta_{26} & \Theta_{27} \\
               \Theta_{31} & \Theta_{32} & \Theta_{33} & \Theta_{34} & \Theta_{35} & \Theta_{36} & \Theta_{37} \\
               \Theta_{41} & \Theta_{42} & \Theta_{43} & \Theta_{44} & \Theta_{45} & \Theta_{46} & \Theta_{47} \\
               \Theta_{51} & \Theta_{52} & \Theta_{53} & \Theta_{54} & \Theta_{55} & \Theta_{56} & \Theta_{57} \\
               \Theta_{61} & \Theta_{62} & \Theta_{63} & \Theta_{64} & \Theta_{65} & \Theta_{66} & \Theta_{67} \\
               \Theta_{71} & \Theta_{72} & \Theta_{73} & \Theta_{74} & \Theta_{75} & \Theta_{76} & \Theta_{77} \\
             \end{array}
           \right),
\end{equation}
Similar to the above scattering coefficients $\Omega_{kj}$,
we can know that $\Theta_{77}$ allows an analytical extension to $\mathbb{C}^{-}$.
Otherwise $\Theta_{kj}(1\leq k,j\leq 6)$ are analytically extendible to $\mathbb{C}^{+}$.
Besides, $\Theta_{k7}, r_{7j}(1\leq k,j\leq 6)$ can be defined on thee real $\lambda$-axis.

Summarizing the above analysis, we have derived two matrix functions $P_{1}$ and $P_{2}$,
which are analytic in $\mathbb{C}^{+}$ and  $\mathbb{C}^{-}$, respectively.
Next we notice that the limit of $P_{2}$ from the right-hand side of the real $\lambda$-line as $P^{-1}$,
and the limit of $P_{1}$ from the left-hand side the real $\lambda$-line as $P^{+}$,
Hence we can get an RHP for the Eq. \eqref{SSE}
\begin{equation}\label{RH-19}
G(x,\lambda)=P^{-}(x,\lambda)P^{+}(x,\lambda),~~\lambda\in\mathbb{R},
\end{equation}
where
%{\small\begin{equation}\label{RH-20}
%G=\left(
%    \begin{array}{cc}
%      I_{6\times 6} & \Omega_{i7}e^{2i\lambda x} \\
%      \Omega_{7i}e^{2i\lambda x} & 1 \\
%    \end{array}
%  \right),~~i=1,2,3,4,5,6,
%\end{equation}}
\begin{align}\label{RH-20}
G=\left(
     \begin{array}{cc}
       I_{6\times6} & \Omega_{j7}e^{2i\lambda x} \\
       \Theta_{7j}e^{2i\lambda x} & 1 \\
     \end{array}
   \right),~~j=1,2,\ldots,6,
\end{align}
and the canonical normalization condition for above RHP \eqref{RH-19} yields
\begin{equation}\label{RH-21}
\left\{ \begin{aligned}
&P_{1}(x,\lambda)\rightarrow\mathbb{I},~~\lambda\in\mathbb{C}^{+}\rightarrow\infty,\\
&P_{2}(x,\lambda)\rightarrow\mathbb{I},~~\lambda\in\mathbb{C}^{-}\rightarrow\infty.
             \end{aligned} \right.
\end{equation}

Next we suppose that the RHP \eqref{RH-19} is irregular.
Here it indicates both $\det(P_{1})$ and $\det(P_{2})$ have certain zero in the analytic domains.
From the definitions of $P_{1}$ and $P_{2}$, we have
\begin{equation}\label{RH-22}
\det(P_{1}(\lambda))=\Omega_{77}(\lambda)~~~~\mbox{and}~~~~\det(P_{2}(\lambda))=\Theta_{77}(\lambda).
\end{equation}
To explain these zero well, we introduce a symmetry relation $Q^{\dag}=-Q$,
%\begin{equation}\label{RH-23}
%Q^{\dag}=-Q,
%\end{equation}
where ``\dag'' means the Hermitian of a matrix.
Hence from the relation \eqref{RHP-12}, we have
\begin{equation}\label{RH-24}
\Psi_{\pm}^{\dag}(\bar{\lambda})=\Psi_{\pm}^{-1}(\lambda).
\end{equation}
The it follows from \eqref{RH-24} that
%\begin{equation}\label{RH-25}
%\Omega^{\dag}(\bar{\lambda})=\Omega^{-1}(\lambda),
%\end{equation}
%which gives the following relations
\begin{align}\label{RH-26}
\Omega^{\dag}(\bar{\lambda})=\Omega^{-1}(\lambda)\rightarrow
\left\{ \begin{aligned}
&\bar{\Omega}_{17}(\lambda)=\Theta_{71}(\lambda),~~\bar{\Omega}_{27}(\lambda)=\Theta_{72}(\lambda),~~\lambda\in\mathbb{R},\\
&\bar{\Omega}_{37}(\lambda)=\Theta_{73}(\lambda),~~\bar{\Omega}_{47}(\lambda)=\Theta_{74}(\lambda),~~\lambda\in\mathbb{R},\\
&\bar{\Omega}_{57}(\lambda)=\Theta_{75}(\lambda),~~\bar{\Omega}_{67}(\lambda)=\Theta_{76}(\lambda),~~\lambda\in\mathbb{R},\\
&\Omega_{77}(\lambda)=\bar{\Theta}_{55}(\bar{\lambda}),~~\lambda\in\mathbb{C}^{+}.
                  \end{aligned} \right.
\end{align}
In addition, the following expression also holds
\begin{equation}\label{RH-27}
P_{1}^{\dag}(\bar{\lambda})=P_{2}(\lambda),~~~~\lambda\in\mathbb{C}^{-}.
\end{equation}
To further settle the RHP \eqref{RH-19}, we must introduce another symmetry relation $\bar{Q}=\sigma Q\sigma$,
%\begin{equation}\label{RH-28}
%\bar{Q}=\sigma Q\sigma,
%\end{equation}
where
\begin{equation}\label{RH-29}
\sigma=\left(
         \begin{array}{ccccccc}
           0 & 1 & 0 & 0 & 0 & 0 & 0 \\
           1 & 0 & 0 & 0 & 0 & 0 & 0 \\
           0 & 0 & 0 & 1 & 0 & 0 & 0 \\
           0 & 0 & 1 & 0 & 0 & 0 & 0 \\
           0 & 0 & 0 & 0 & 0 & 1 & 0 \\
           0 & 0 & 0 & 0 & 1 & 0 & 0 \\
           0 & 0 & 0 & 0 & 0 & 0 & 1 \\
         \end{array}
       \right).
\end{equation}
It follows from \eqref{RH-29} that
\begin{equation}\label{RH-30}
\sigma\bar{\Psi}_{\pm}(-\bar{\lambda})\sigma=\Psi_{\pm}(\lambda),\rightarrow\sigma \bar{\Omega}(-\bar{\lambda})\sigma=\Omega(\lambda).
\end{equation}
%which yields
%\begin{equation}\label{RH-31}
%\sigma \bar{\Omega}(-\bar{\lambda})\sigma=\Omega(\lambda).
%\end{equation}
Apparently, the above relation means that
\begin{equation}\label{RH-32}
\left\{ \begin{aligned}
&\Omega_{17}(\lambda)=\bar{\Omega}_{27}(-\lambda),~~\Omega_{37}(\lambda)=\bar{\Omega}_{47}(-\lambda),~~\Omega_{57}(\lambda)=\bar{\Omega}_{67}(-\lambda),~~\lambda\in\mathbb{R},\\
&\Omega_{77}(\lambda)=\bar{\Omega}_{77}(-\bar{\lambda}),~~\lambda\in\mathbb{C}^{+}.
                     \end{aligned} \right.
\end{equation}
In addition, from the relation \eqref{RH-30}, we have
\begin{equation}\label{RH-33}
\sigma \bar{P}_{1}(-\bar{\lambda})\sigma=P_{1}(\lambda),~~\lambda\in\mathbb{C}^{+}.
\end{equation}
Thus, we know that if $\lambda_{j}$ is a zero of $\det(P_{1})$, then $\hat{\lambda}_{j}=\bar{\lambda}_{j}$ is a zero of $\det(P_{2})$.
Furthermore, in terms of \eqref{RH-32}, we find that $-\bar{\lambda}_{j}$ is also a zero of  $\det(P_{1})$.
Next we consider the zeros of  $\det(P_{1})$ in the following two cases.\\
\textbf{Case (I):} we suppose that  $\det(P_{1})$ admits a total number of $2N$ zeros $\lambda_{j}(1\leq j\leq 2N)$
satisfying $\lambda_{N+j}=-\bar{\lambda}_{j}(1\leq j\leq N)$, which are all in $\mathbb{C}^{+}$.
Likewise, $\det(P_{2})$ admits $2N$ zeros $\hat{\lambda}_{j}$ $(1\leq j\leq 2N)$ satisfy $\hat{\lambda}_{j}=\bar{\lambda}_{j}$,
which all lie in $\mathbb{C}^{-}$.\\
\textbf{Case (II):} $\det(P_{1})$ admits also $N$ simple zeros $\lambda_{j} (1\leq j\leq N)$ in $\mathbb{C}^{+}$,
where all $\lambda_{j}$ are pure imaginary.
$\det(P_{2})$ admits $N$ zeros $\hat{\lambda}_{j}$ in $\mathbb{C}^{-}$, where $\hat{\lambda}_{j}=\bar{\lambda}_{j}$.
Under these assumptions, ker($P_{1}(\lambda_{j})$ and $P_{2}(\hat{\lambda}_{j})$ are one-dimensional, and they are spanned by
\begin{equation}\label{RH-34}
P_{1}(\lambda_{j})v_{j}=0,~~\hat{v}_{j}P_{2}(\hat{\lambda}_{j})=0,~~1\leq j\leq2N,
\end{equation}
respectively. For the first type of zeros, from the relations \eqref{RH-27} and \eqref{RH-34}, we have $\hat{v}_{j}=v_{j}^{\dag},1\leq j\leq 2N$.
%\begin{equation}\label{RH-35}
%\hat{v}_{j}=v_{j}^{\dag},~~~~1\leq j\leq 2N.
%\end{equation}
Similarly, form \eqref{RH-33}, we can obtain the following relation $v_{j}=\sigma \bar{v}_{j-N},N+1\leq j\leq2N$.
%\begin{equation}\label{RH-36}
%v_{j}=\sigma \bar{v}_{j-N},~~~~N+1\leq j\leq2N.
%\end{equation}
Next we should construct the vectors $v_{j} (1\leq j\leq N)$.
To this end, we set the x-derivative of $P_{1}(\lambda_{j})v_{j}=0$.
The using the first expression in \eqref{RHP-2}, we can obtain $v_{j}=e^{i\lambda_{j}\sigma_{3}x}v_{j,0},1\leq j\leq N$,
%\begin{equation}\label{RH-37}
%v_{j}=e^{i\lambda_{j}\sigma_{3}x}v_{j,0},~~~~1\leq j\leq N,
%\end{equation}
in which $v_{j0}$ is independent of x.
Thus utilizing the above results, the two vectors $\hat{v}_{j}$ and $v_{j}$ can be obtained explicitly.
To construct the multi-soliton solutions for the Eq. \eqref{SSE},
we should choose $G=\mathbb{I}$.
Therefore, the solution for the particular RHP reads
\begin{equation}\label{RH-38}
\left\{ \begin{aligned}
&P_{1}(\lambda)=\mathbb{I}-\sum_{k=1}^{2N}\sum_{j=1}^{2N}\frac{v_{k}v_{j}(M^{-1})_{kj}}{\lambda-\hat{\lambda}_{j}},\\
&P_{2}(\lambda)=\mathbb{I}+\sum_{k=1}^{2N}\sum_{j=1}^{2N}\frac{v_{k}v_{j}(M^{-1})_{kj}}{\lambda-\hat{\lambda}_{k}},
                          \end{aligned} \right.
\end{equation}
where $M=(M_{kj})_{2N\times 2N}$ is a matrix whose entries yield $M_{kj}=\frac{\hat{v}_{k}k_{j}}{\lambda_{j}-\bar{\lambda}_{k}}$.
%\begin{equation}\label{RH-39}
%M_{kj}=\frac{\hat{v}_{k}k_{j}}{\lambda_{j}-\bar{\lambda}_{k}}.
%\end{equation}
For second type of zeros, the vectors $v_{j}, \hat{v}_{j} (1\leq j\leq N)$ are defined by
\begin{equation}\label{RH-40}
\hat{v}_{j}=v_{j}^{\dag},~~v_{j}=e^{i\lambda_{j}\sigma_{3}x}v_{j,0}.
\end{equation}
Utilizing these vectors, the RHP \eqref{RH-19} in this case can be also solved
\begin{equation}\label{RH-41}
\left\{ \begin{aligned}
&P_{1}(\lambda)=\mathbb{I}-\sum_{k=1}^{N}\sum_{j=1}^{N}\frac{v_{k}v_{j}(M^{-1})_{kj}}{\lambda-\hat{\lambda}_{j}},\\
&P_{2}(\lambda)=\mathbb{I}+\sum_{k=1}^{N}\sum_{j=1}^{N}\frac{v_{k}v_{j}(M^{-1})_{kj}}{\lambda-\hat{\lambda}_{k}},
  \end{aligned} \right.
\end{equation}
where $M=(M_{kj})_{N\times N}$ is a matrix whose entries yield $M_{kj}=\frac{\hat{v}_{k}k_{j}}{\lambda_{j}-\bar{\lambda}_{k}}$.
%\begin{equation}\label{RH-39}
%M_{kj}=\frac{\hat{v}_{k}k_{j}}{\lambda_{j}-\bar{\lambda}_{k}}.
%\end{equation}
In terms of $P_{1}$ in \eqref{RH-38}, we can recover the potentials $u_{1}$, $u_{2}$, $u_{3}$.
Actually, we can expand $P_{1}(\lambda)$ as
\begin{equation}\label{RH-40}
P_{1}(\lambda)=\mathbb{I}+\frac{P_{1}^{(1)}}{\lambda}+\frac{P_{1}^{(2)}}{\lambda^2}+\left(\frac{1}{\lambda^3}\right),~~\lambda\rightarrow\infty.
\end{equation}
Then inserting \eqref{RH-40} into the first expression in \eqref{RHP-2}, and collecting O(1) terms yields $Q=i\left[P_{1}^{(1)},\sigma_{3}\right]$,
%\begin{equation}\label{RH-41}
%Q=i\left[P_{1}^{(1)},\sigma_{3}\right],
%\end{equation}
which means that $u_{1}$, $u_{2}$, $u_{3}$ can be expressed as
\begin{align}\label{RH-42}
&u_{1}=-2i\left(P_{1}^{(1)}\right)_{17},~~u_{2}=-2i\left(P_{1}^{(1)}\right)_{37},~~
u_{3}=-2i\left(P_{1}^{(1)}\right)_{57},
\end{align}
where $(P_{1}^{(1)})_{k7}$ is the (k,7)-entry of the function $P_{1}^{(1)}$.
Here from \eqref{RH-38} and \eqref{RH-41}, the matrix functions can be rewritten as
\begin{equation}\label{RH-43}
\left\{ \begin{aligned}
&P_{1}^{(1)}=-\sum_{k=1}^{2N}\sum_{j}^{2N}v_{k}v_{j}\left(M^{-1}\right)_{kj},\\
&P_{1}^{(1)}=-\sum_{k=1}^{N}\sum_{j}^{N}v_{k}v_{j}\left(M^{-1}\right)_{kj}.
      \end{aligned} \right.
\end{equation}

\section{Multi-soliton solutions}
To obtain soliton solutions for the Eq. \eqref{SSE}, we must consider the t-evolutions of the scattering data.
From the second expression in \eqref{RHP-2} and \eqref{RHP-8}, we get $\Omega_{t}=4i\lambda^{3}[\sigma_{3},\Omega]$,
%\begin{equation}\label{SS-1}
% \Omega_{t}=4i\lambda^{3}[\sigma_{3},\Omega],
%\end{equation}
which gives the following results
\begin{equation}\label{SS-2}
\Omega_{17,t}=8i\lambda^3 \Omega_{15},~~\Omega_{37,t}=8i\lambda^3\Omega_{37},~~~~\Omega_{57,t}=8i\lambda^3\Omega_{57}.
\end{equation}
Additionally, utilizing the second expression in \eqref{RHP-2}, we can obtain that $v_{j,t}=4i\lambda_{j}^3\sigma_{3}v_{j}$.
%\begin{equation}\label{SS-3}
%v_{j,t}=4i\lambda_{j}^3\sigma_{3}v_{j}.
%\end{equation}
Thus, for the first type of zeros, we have
\begin{equation}\label{SS-4}
v_{j}=
\left\{ \begin{aligned}
&e^{\theta_{j}\sigma}v_{j,0},~~~~~~~~~~~~1\leq j\leq N,\\
&\sigma e^{\bar{\theta}_{j-N}\sigma_{3}}\bar{v}_{j-N,0},~~N+1\leq j\leq2N,
            \end{aligned} \right.
\end{equation}
and
\begin{equation}\label{SS-5}
\hat{v}_{j}=
\left\{ \begin{aligned}
&v_{j,0}^{\dag}e^{\bar{\theta}_{j}\sigma},~~~~~~~~~~~~1\leq j\leq N,\\
&v_{j-N,0}^{T}e^{\bar{\theta}_{j-N}\sigma_{3}}\sigma,~~N+1\leq j\leq2N,
           \end{aligned} \right.
\end{equation}
where $\theta_{j}=i\lambda_{j}x+4i\lambda_{j}^3t$, and $v_{j,0}$ are constant vectors.
Furthermore, for the second type of zeros, we have
\begin{equation}\label{SS-6}
v_{j}=e^{\theta_{j}\sigma_{3}}v_{j,0},~~\hat{v}_{j}=v_{j,0}^{\dag}e^{\bar{\theta}_{j}\sigma_{3}},~~1\leq j\leq N,
\end{equation}
where $\theta_{j}=i\lambda_{j}x+4i\lambda_{j}^3t$ with $\lambda_{j}$ being imaginary, and $v_{j,0}$ are constant vectors.
For the fist type of zeros of $\det(P_{1})$,
we choose $v_{j,0}=(\alpha_{j},\beta_{j},\gamma_{j},\mu_{j},\rho_{j},\delta_{j},1)^{T}$
to be complex constant vectors.
Then using Eqs.\eqref{RH-42},\eqref{RH-43}, \eqref{SS-4} and \eqref{SS-5},
the $N$-soliton solutions of the Eq. \eqref{SSE} can be obtained as follows
\begin{equation}\label{SS-7}
\left\{ \begin{aligned}
u_{1}=&2i\sum_{k=1}^{N}\sum_{j=1}^{N}\alpha_{k}e^{\theta_{k}-\bar{\theta}_{j}}(M^{-1})_{kj}
+2i\sum_{k=1}^{N}\sum_{j=N+1}^{2N}\alpha_{k}e^{\theta_{k}-\theta_{j-N}}(M^{-1})_{kj}\\
&+2i\sum_{k=N+1}^{2N}\sum_{j=1}^{N}\bar{\beta}_{k-N}e^{\bar{\theta}_{k-N}-\bar{\theta}_{j}}(M^{-1})_{kj}
+2i\sum_{k=N+1}^{2N}\sum_{j=N+1}^{2N}\bar{\beta}_{k-N}e^{\bar{\theta}_{k-N}-\theta_{j-N}}(M^{-1})_{kj},\\
u_{2}=&2i\sum_{k=1}^{N}\sum_{j=1}^{N}\gamma_{k}e^{\theta_{k}-\bar{\theta}_{j}}(M^{-1})_{kj}
+2i\sum_{k=1}^{N}\sum_{j=N+1}^{2N}\gamma_{k}e^{\theta_{k}-\theta_{j-N}}(M^{-1})_{kj}\\
&+2i\sum_{k=N+1}^{2N}\sum_{j=1}^{N}\bar{\mu}_{k-N}e^{\bar{\theta}_{k-N}-\bar{\theta}_{j}}(M^{-1})_{kj}
+2i\sum_{k=N+1}^{2N}\sum_{j=N+1}^{2N}\bar{\mu}_{k-N}e^{\bar{\theta}_{k-N}-\theta_{j-N}}(M^{-1})_{kj},\\
u_{3}=&2i\sum_{k=1}^{N}\sum_{j=1}^{N}\rho_{k}e^{\theta_{k}-\bar{\theta}_{j}}(M^{-1})_{kj}
+2i\sum_{k=1}^{N}\sum_{j=N+1}^{2N}\rho_{k}e^{\theta_{k}-\theta_{j-N}}(M^{-1})_{kj}\\
&+2i\sum_{k=N+1}^{2N}\sum_{j=1}^{N}\bar{\delta}_{k-N}e^{\bar{\theta}_{k-N}-\bar{\theta}_{j}}(M^{-1})_{kj}
+2i\sum_{k=N+1}^{2N}\sum_{j=N+1}^{2N}\bar{\delta}_{k-N}e^{\bar{\theta}_{k-N}-\theta_{j-N}}(M^{-1})_{kj},\\
               \end{aligned} \right.
\end{equation}
where $M=(M_{kj})_{2N\times 2N}$ is defined by
\begin{align}\label{SS-8}
M_{kj}=
&\left\{ \begin{aligned}
&\frac{\Delta_{1}e^{\bar{\theta}_{k}+\theta_{j}}+
e^{\bar{\theta}_{k}+\theta_{j}}}{\lambda_{j}-\bar{\lambda}_{k}},~1\leq k,j\leq N,\\
&\frac{\Delta_{2}e^{\bar{\theta}_{k}+\bar{\theta}_{j-N}}+
e^{-\bar{\theta}_{k}-\bar{\theta}_{j-N}}}{-\bar{\lambda}_{j-N}-\bar{\lambda}_{k}},~~1\leq k\leq N,~~N+1\leq j\leq2N,\\
&\frac{\Delta_{3}e^{\theta_{j}+\theta_{k-N}}+
e^{-\theta_{j}-\theta_{k-N}}}{\lambda_{k-N}+\lambda_{j}},~~1\leq j\leq N,~N+1\leq k\leq2N,\\
&\frac{\Delta_{4}e^{\bar{\theta}_{j-N}+\theta_{k-N}}+
e^{-\bar{\theta}_{j-N}+\theta_{k-N}}}{\lambda_{k-N}-\bar{\lambda}_{j-N}},~~N+1\leq k,j\leq 2N,\\
                   \end{aligned} \right.
\end{align}
with
\begin{equation}\label{SS-9}
\left\{ \begin{aligned}
&\Delta_{1}=\bar{\alpha}_{k}\alpha_{j}+\bar{\beta}_{k}\beta_{j}+\bar{\gamma}_{k}\gamma_{j}+\bar{\mu}_{k}\mu_{j}+\bar{\rho}_{k}\rho_{j}+\bar{\delta}_{k}\delta_{j},\\
&\Delta_{2}=\bar{\alpha}_{k}\bar{\beta}_{j-N}+\bar{\beta}_{k}\bar{\alpha}_{j-N}+\bar{\gamma}_{k}\bar{\mu}_{j-N}+\bar{\mu}_{k}\bar{\gamma}_{j-N}+\bar{\rho}_{j-N}\bar{\delta}_{k}+\bar{\rho}_{k}\bar{\delta}_{j-N},\\
&\Delta_{3}=\bar{\alpha}_{j}\bar{\beta}_{k-N}+\bar{\beta}_{j}\bar{\alpha}_{k-N}+\bar{\gamma}_{j}\bar{\mu}_{k-N}
+\bar{\mu}_{j}\bar{\gamma}_{k-N}+\bar{\rho}_{k-N}\bar{\delta}_{j}+\bar{\rho}_{j}\bar{\delta}_{k-N},\\
&\Delta_{4}=\bar{\alpha}_{j-N}\alpha_{k-N}+\bar{\beta}_{j-N}\beta_{k-N}+\bar{\gamma}_{j-N}\gamma_{k-N}+\bar{\mu}_{j-N}\mu_{k-N}+\bar{\rho}_{j-N}\rho_{k-N}+\bar{\delta}_{k}\delta_{j}.
                       \end{aligned} \right.
\end{equation}

To show the one-soliton solution explicitly, we should choose the appropriate parameters in \eqref{SS-7}-\eqref{SS-9}.
In the following, we choose $\bar{\alpha}_{1}=\beta_{1}$, $\bar{\gamma}_{1}=\mu_{1}$ and $\bar{\delta}_{1}=\rho_{1}$ in \eqref{SS-7}-\eqref{SS-9}.
Then by taking $\lambda_{1}=\xi_{1}+i\eta_{1}$, a new breather solution for the Eq. \eqref{SSE} can be constructed
\begin{equation}\label{SS-10}
\left\{ \begin{aligned}
&u_{1}=-\frac{2\sqrt{2}\alpha_{1}\xi_{1}\eta_{1}}{\sqrt{|\alpha_{1}|^{2}+|\gamma_{1}|^{2}+|\rho_{1}|^{2}}}\frac{\xi_{1}\cosh(X_{1})\cos(Y_{1})+\eta_{1}\sinh(X_{1})\sin(Y_{1})}{\xi_{1}^{2}\cosh(X_{1})^{2}+\eta_{1}\sin(Y_{1})^{2}},\\
&u_{2}=-\frac{2\sqrt{2}\gamma_{1}\xi_{1}\eta_{1}}{\sqrt{|\alpha_{1}|^{2}+|\gamma_{1}|^{2}+|\rho_{1}|^{2}}}\frac{\xi_{1}\cosh(X_{1})\cos(Y_{1})+\eta_{1}\sinh(X_{1})\sin(Y_{1})}{\xi_{1}^{2}\cosh(X_{1})^{2}+\eta_{1}\sin(Y_{1})^{2}},\\
&u_{3}=-\frac{2\sqrt{2}\rho_{1}\xi_{1}\eta_{1}}{\sqrt{|\alpha_{1}|^{2}+|\gamma_{1}|^{2}+|\rho_{1}|^{2}}}\frac{\xi_{1}\cosh(X_{1})\cos(Y_{1})+\eta_{1}\sinh(X_{1})\sin(Y_{1})}{\xi_{1}^{2}\cosh(X_{1})^{2}+\eta_{1}\sin(Y_{1})^{2}},\\
                          \end{aligned} \right.
\end{equation}
where
\begin{align}\label{SS-11}
&X_{1}=-2\eta_{1}\left(x+4\left(3\xi_{1}^{2}-\eta_{1}^{2}\right)t\right)+\ln\sqrt{2\left(|\alpha_{1}|^{2}+|\gamma_{1}|^{2}+|\rho_{1}|^{2}\right)},\notag\\
&Y_{1}=2\xi_{1}\left(x+4\left(\xi_{1}^{2}-3\eta_{1}^2\right)t\right).
\end{align}
Fig.1 displays the dynamics of periodic breather-type solution \eqref{SS-10}.

We next discuss the case for $N=2$ in Eqs.\eqref{SS-7}-\eqref{SS-9}.
As seen in Fig.2, the two solitons pass
through each other, and their polarizations do not change.
After collision when this right soliton passes to the left, its power has diminished
dramatically. Particularly, Figs.2(b) and 2(c) also display that the components of $u_{2},u_{3}$ even vanishes after collision.
The another collision is displayed in Figs.2(d)-2(f)).
This collision is a little similar to that in Fig.2(a) and is another example of
soliton interactions. The shapes of the solitons can be changed after collision.
All the phenomena indicate that there is a lot of
energy transfer has taken place between these two solitons during the collision.
Due to the limit of length, the other corresponding figures are omitted here.

%\noindent
%{\rotatebox{0}{\includegraphics[width=4.5cm,height=3.9cm,angle=0]{3-1.eps}}}
%~~~
%{\rotatebox{0}{\includegraphics[width=4.5cm,height=3.9cm,angle=0]{3-2.eps}}}
%~~~
%{\rotatebox{0}{\includegraphics[width=4.5cm,height=3.9cm,angle=0]{3-3.eps}}}
%
%$\qquad\qquad\textbf{(a)}\qquad\qquad\qquad\qquad\qquad\qquad\textbf{(b)}
%\qquad\qquad\qquad\qquad\qquad\qquad\textbf{(c)}$\\
%
%\noindent { \small \textbf{Figure 3.} (Color online) Two-soliton solutions \eqref{SS-10} (N=2)
%with parameters: $\lambda_{1}=0.5+0.5i, \lambda_{2}=0.4+0.6i,  \alpha_{1}=\beta_{1}=\gamma_{1}=\mu_{1}=\rho_{1}=\alpha_{2}=1, \gamma_{2}=2, \beta_{2}=\mu_{2}=\rho_{2}=0$.\\}

In the following, we assume that $\det(P_{1})$ admits $N$ simple zeros $\lambda_{j}$ in $\mathbb{C}^{+}$.
Taking $v_{j,0}=(\alpha_{j},\bar{\alpha}_{j},\gamma_{j},\bar{\gamma}_{j},\rho_{j},\bar{\rho}_{j},1)$ to be complex vectors.
Then by utilizing \eqref{RH-42}, \eqref{RH-42}, and \eqref{SS-5}, we can obtain another kind of $N$-soliton solutions as follows
\begin{equation}\label{SS-12}
\left\{ \begin{aligned}
&u_{1}(x,t)=2i\sum_{k=1}^{N}\sum_{j=1}^{N}\alpha_{k}e^{\theta_{k}-\bar{\theta}_{j}}(M^{-1})_{kj},\\
&u_{2}(x,t)=2i\sum_{k=1}^{N}\sum_{j=1}^{N}\gamma_{k}e^{\theta_{k}-\bar{\theta}_{j}}(M^{-1})_{kj},\\
&u_{3}(x,t)=2i\sum_{k=1}^{N}\sum_{j=1}^{N}\rho_{k}e^{\theta_{k}-\bar{\theta}_{j}}(M^{-1})_{kj},
                          \end{aligned} \right.
\end{equation}
where $M=(M_{kj})_{N\times N} (1\leq k,j\leq N)$ is defined by
\begin{align}\label{SS-13}
M_{kj}=
\frac{(\alpha_{k}\bar{\alpha}_{j}+\bar{\alpha}_{k}\alpha_{j}+\gamma_{k}\bar{\gamma}_{j}
+\bar{\gamma}_{k}\gamma_{j}+\bar{\rho}_{k}\rho_{j}+\rho_{k}\bar{\rho}_{j})e^{\bar{\theta}_{k}+\theta_{j}}+e^{-\bar{\theta}_{k}-\theta_{j}}}
{\lambda_{j}-\bar{\lambda}_{k}}.
\end{align}
Taking $N=1$, the $N$-soliton solutions \eqref{SS-12} yields a single soliton solution for the Eq. \eqref{SSE}
\begin{equation}\label{SS-14}
\left\{ \begin{aligned}
&u_{1}=-\frac{\sqrt{2}\alpha_{1}\eta_{1}}{\sqrt{|\alpha_{1}|^{2}+|\gamma_{1}|^{2}+|\rho_{1}|^{2}}}\mbox{sech}\left(2\eta_{1}x-8\eta_{1}^{3}t+\ln\left(\sqrt{|\alpha_{1}|^{2}+|\gamma_{1}|^{2}+|\rho_{1}|^{2}}\right)\right),\\
&u_{2}=-\frac{\sqrt{2}\gamma_{1}\eta_{1}}{\sqrt{|\alpha_{1}|^{2}+|\gamma_{1}|^{2}+|\rho_{1}|^{2}}}\mbox{sech}\left(2\eta_{1}x-8\eta_{1}^{3}t+\ln\left(\sqrt{|\alpha_{1}|^{2}+|\gamma_{1}|^{2}+|\rho_{1}|^{2}}\right)\right),\\
&u_{3}=-\frac{\sqrt{2}\rho_{1}\eta_{1}}{\sqrt{|\alpha_{1}|^{2}+|\gamma_{1}|^{2}+|\rho_{1}|^{2}}}\mbox{sech}\left(2\eta_{1}x-8\eta_{1}^{3}t+\ln\left(\sqrt{|\alpha_{1}|^{2}+|\gamma_{1}|^{2}+|\rho_{1}|^{2}}\right)\right),
                             \end{aligned} \right.
\end{equation}
Fig.3(a) and Fig.3(b) present the bright-soliton and dark-soliton solution, respectively, by
seeking the suitable parameters, which is useful for
understanding the dynamical behaviors of the soliton solutions.

Next, taking $N=2$, Eq.\eqref{SS-12} can be reduced to a two-bell soliton solution of the Eq. \eqref{SSE} given by
\begin{equation}\label{Two-2}
\left\{ \begin{aligned}
  &u_{1}(x,t)=2i\alpha_{1}\exp\left(\theta_{1}-\bar{\theta}_{1}\right)\left(M^{-1}\right)_{11}
  +2i\alpha_{1}\exp\left(\theta_{1}-\bar{\theta}_{2}\right)\left(M^{-1}\right)_{12}\\
  &~~~~~~~~~~+2i\alpha_{2}\exp\left(\theta_{2}-\bar{\theta}_{1}\right)\left(M^{-1}\right)_{21}
  +2i\alpha_{2}\exp\left(\theta_{2}-\bar{\theta}_{2}\right)\left(M^{-1}\right)_{22},\\
  &u_{2}(x,t)=2i\gamma_{1}\exp\left(\theta_{1}-\bar{\theta}_{1}\right)\left(M^{-1}\right)_{11}
  +2i\gamma_{1}\exp\left(\theta_{1}-\bar{\theta}_{2}\right)\left(M^{-1}\right)_{12}\\
  &~~~~~~~~~~+2i\gamma_{2}\exp\left(\theta_{2}-\bar{\theta}_{1}\right)\left(M^{-1}\right)_{21}
  +2i\gamma_{2}\exp\left(\theta_{2}-\bar{\theta}_{2}\right)\left(M^{-1}\right)_{22},\\
  &u_{3}(x,t)=2i\rho_{1}\exp\left(\theta_{1}-\bar{\theta}_{1}\right)\left(M^{-1}\right)_{11}
  +2i\rho_{1}\exp\left(\theta_{1}-\bar{\theta}_{2}\right)\left(M^{-1}\right)_{12}\\
  &~~~~~~~~~~+2i\rho_{2}\exp\left(\theta_{2}-\bar{\theta}_{1}\right)\left(M^{-1}\right)_{21}
  +2i\rho_{2}\exp\left(\theta_{2}-\bar{\theta}_{2}\right)\left(M^{-1}\right)_{22},
        \end{aligned} \right.
\end{equation}
where $M=(T_{jm})_{2\times2}$ with
\begin{equation}\label{Two-3}
\left\{ \begin{aligned}
  &T_{11}=\frac{2\left(|\alpha_{1}|^{2}+|\gamma_{1}|^{2}+|\rho_{1}|^2\right)\exp\left(\theta_{1}+\bar{\theta}_{1}\right)
  +\exp\left(-\theta_{1}-\bar{\theta}_{1}\right)}{\lambda_{1}-\bar{\lambda}_{1}},\\
  &T_{12}=\frac{\nabla_{1}\exp\left(\bar{\theta}_{1}+\theta_{2}\right)
  +\exp\left(-\bar{\theta}_{1}-\theta_{2}\right)}
  {\lambda_{2}-\bar{\lambda}_{1}},\\
  &T_{21}=\frac{\nabla_{2}\exp\left(\theta_{1}+\bar{\theta}_{2}\right)
  +\exp\left(-\theta_{1}-\bar{\theta}_{2}\right)}
  {\lambda_{1}-\bar{\lambda}_{2}},\\
  &T_{22}=\frac{2\left(|\alpha_{1}|^{2}+|\gamma_{1}|^{2}+|\rho_{1}|^2\right)\exp\left(\theta_{2}+\bar{\theta}_{2}\right)
  +\exp\left(-\theta_{2}-\bar{\theta}_{2}\right)}{\lambda_{2}-\bar{\lambda}_{2}},\\
  &\nabla_{1}=\left(\bar{\alpha}_{1}\alpha_{2}+\alpha_{1}\bar{\alpha}_{2}
  +\bar{\gamma}_{1}\gamma_{2}+\gamma_{1}\bar{\gamma}_{2}+\bar{\rho}_{1}\rho_{2}+\rho_{1}\bar{\rho}_{2}\right),\\
  &\nabla_{2}=\left(\bar{\alpha}_{1}\alpha_{2}+\alpha_{1}\bar{\alpha}_{2}
  +\bar{\gamma}_{1}\gamma_{2}+\gamma_{1}\bar{\gamma}_{2}+\bar{\rho}_{1}\rho_{2}+\rho_{1}\bar{\rho}_{2}\right),
          \end{aligned} \right.
\end{equation}
and $\theta_{m}=i\lambda_{m}x+4i\lambda_{m}^3t(m=1,2)$.
The two-bell soliton interactions given by Eq.\eqref{SS-12} with $N=2$ are shown in Fig.4.

\noindent
{\rotatebox{0}{\includegraphics[width=4.8cm,height=2.8cm,angle=0]{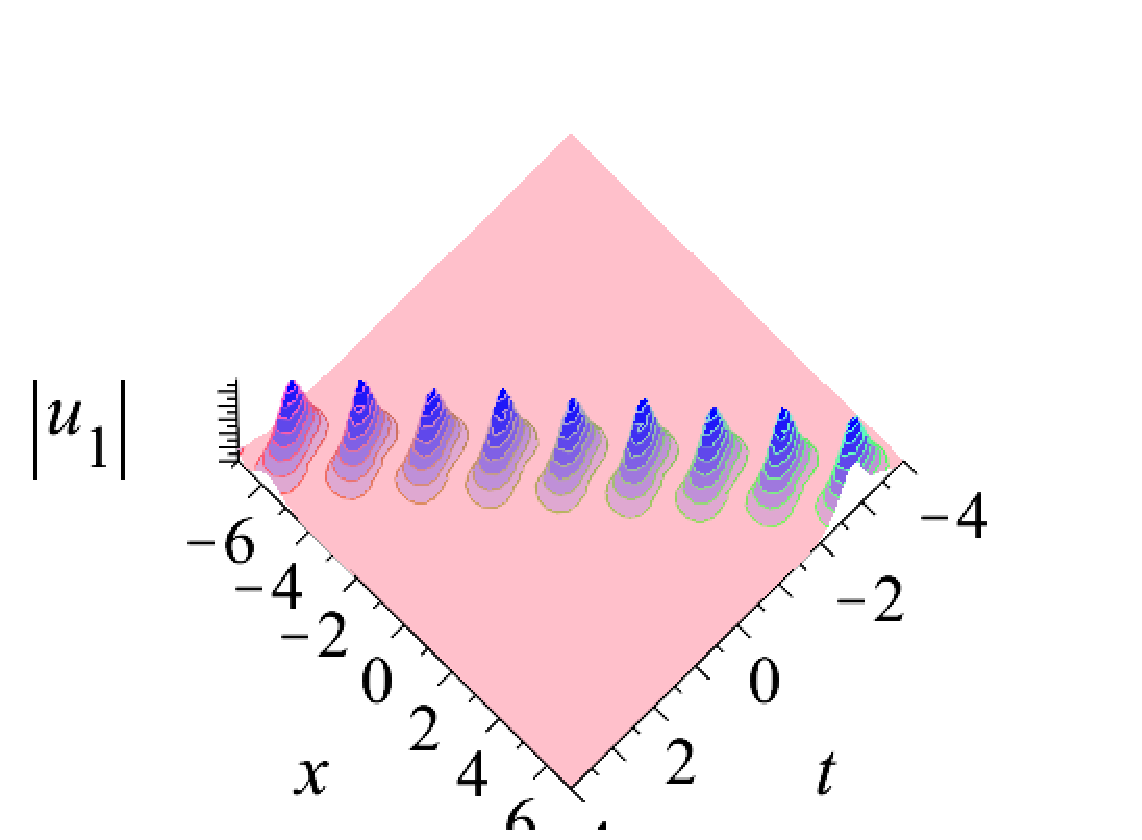}}}
{\rotatebox{0}{\includegraphics[width=4.8cm,height=2.8cm,angle=0]{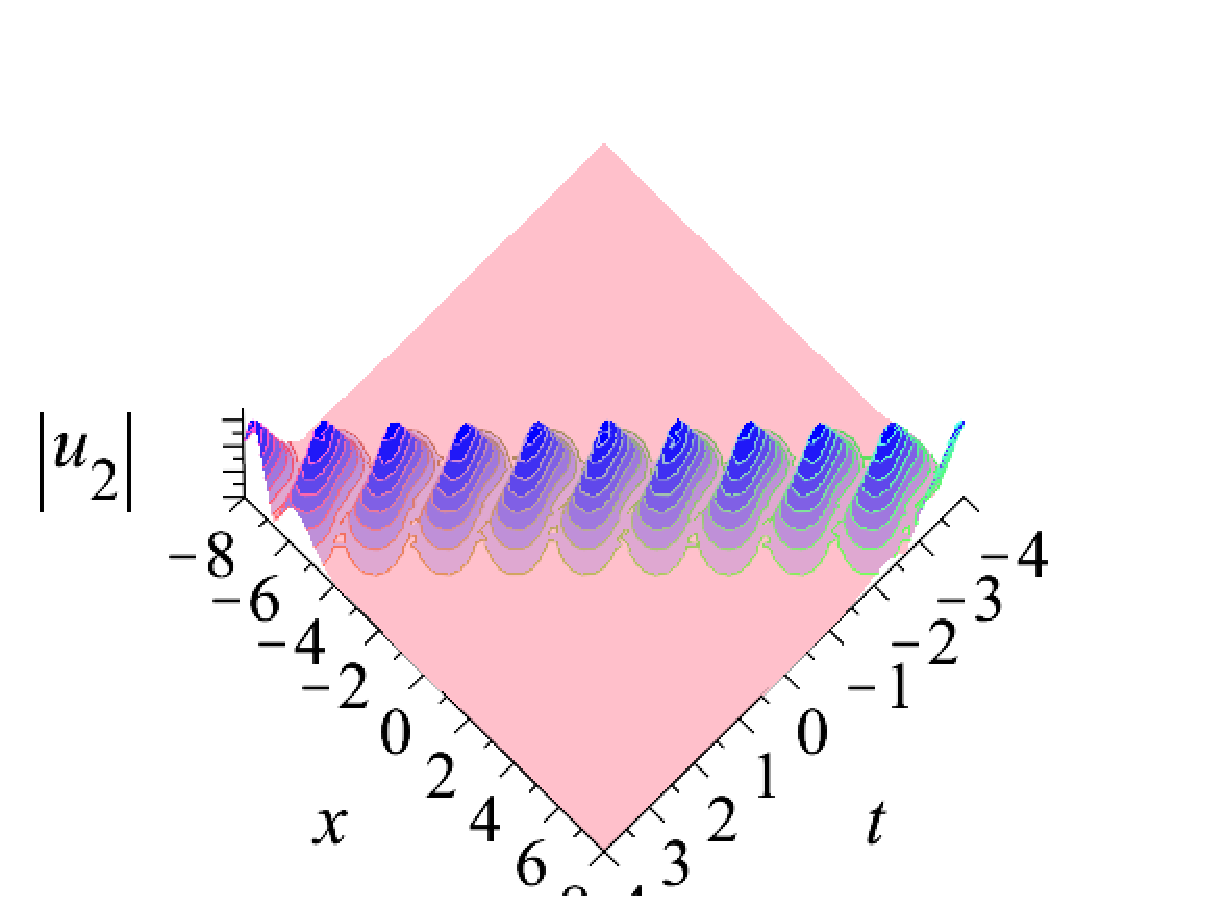}}}
{\rotatebox{0}{\includegraphics[width=4.8cm,height=2.8cm,angle=0]{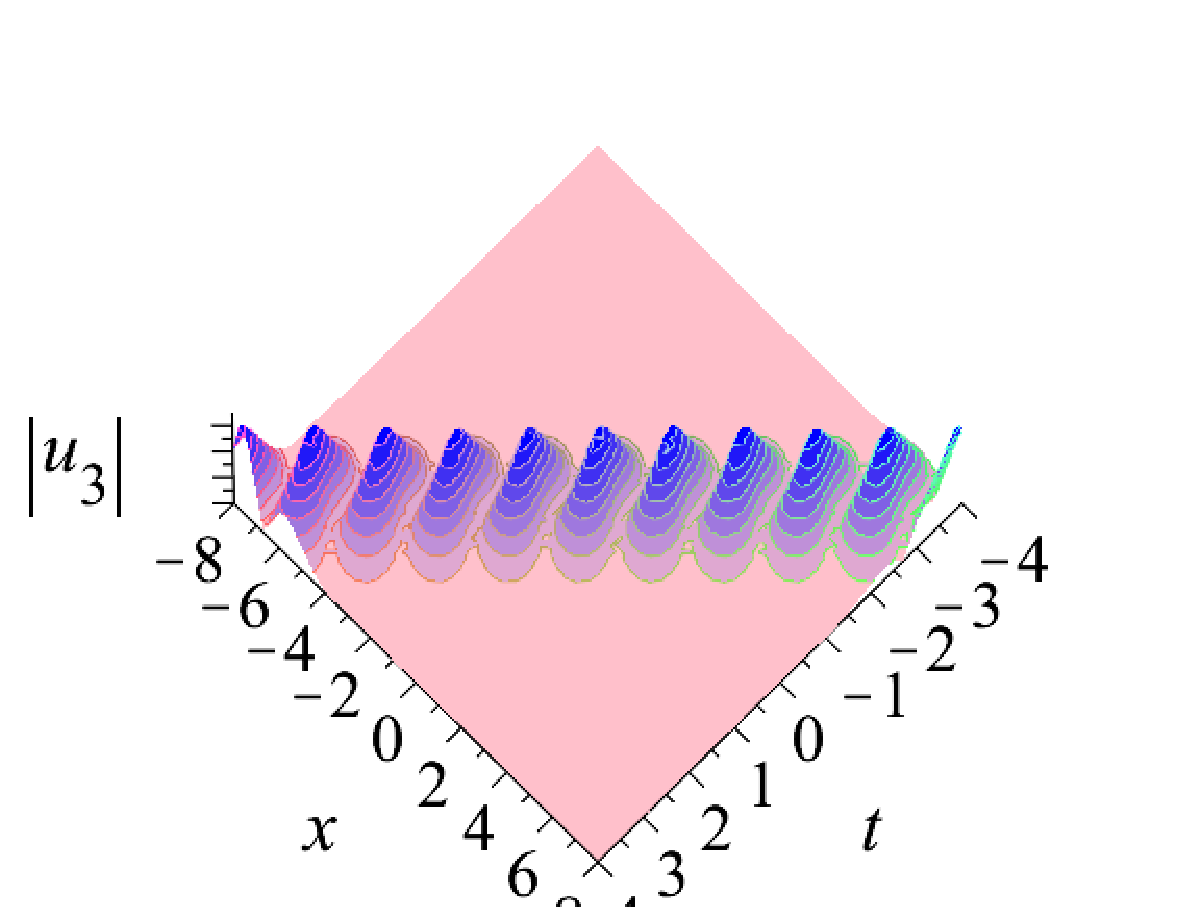}}}

$\qquad\qquad\textbf{(a)}\qquad\qquad\qquad\qquad\qquad\textbf{(b)}
\qquad\qquad\qquad\qquad\qquad\qquad\textbf{(c)}$\\

\noindent { \small \textbf{Figure 1.} (Color online) Breather wave via solutions \eqref{SS-10} $(|u_{1}|$, $|u_{2}|$, $|u_{3}|)$
with parameters: $\alpha_{1}=i/\sqrt{3}, \gamma_{1}=\sqrt{2}i/\sqrt{3}, \rho_{1}=\sqrt{2}i/\sqrt{3}, \lambda_{1}=0.5+0.5i$.\\}

\noindent
{\rotatebox{0}{\includegraphics[width=4.8cm,height=2.8cm,angle=0]{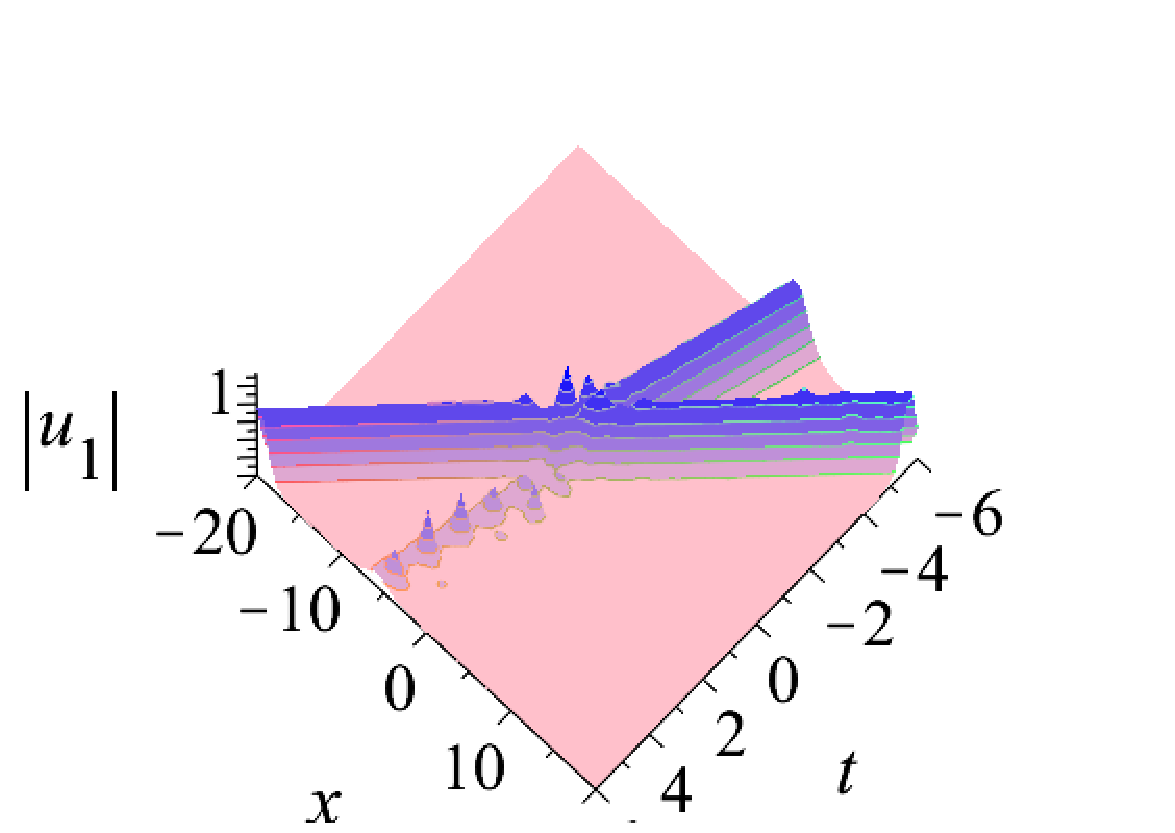}}}
{\rotatebox{0}{\includegraphics[width=4.8cm,height=2.8cm,angle=0]{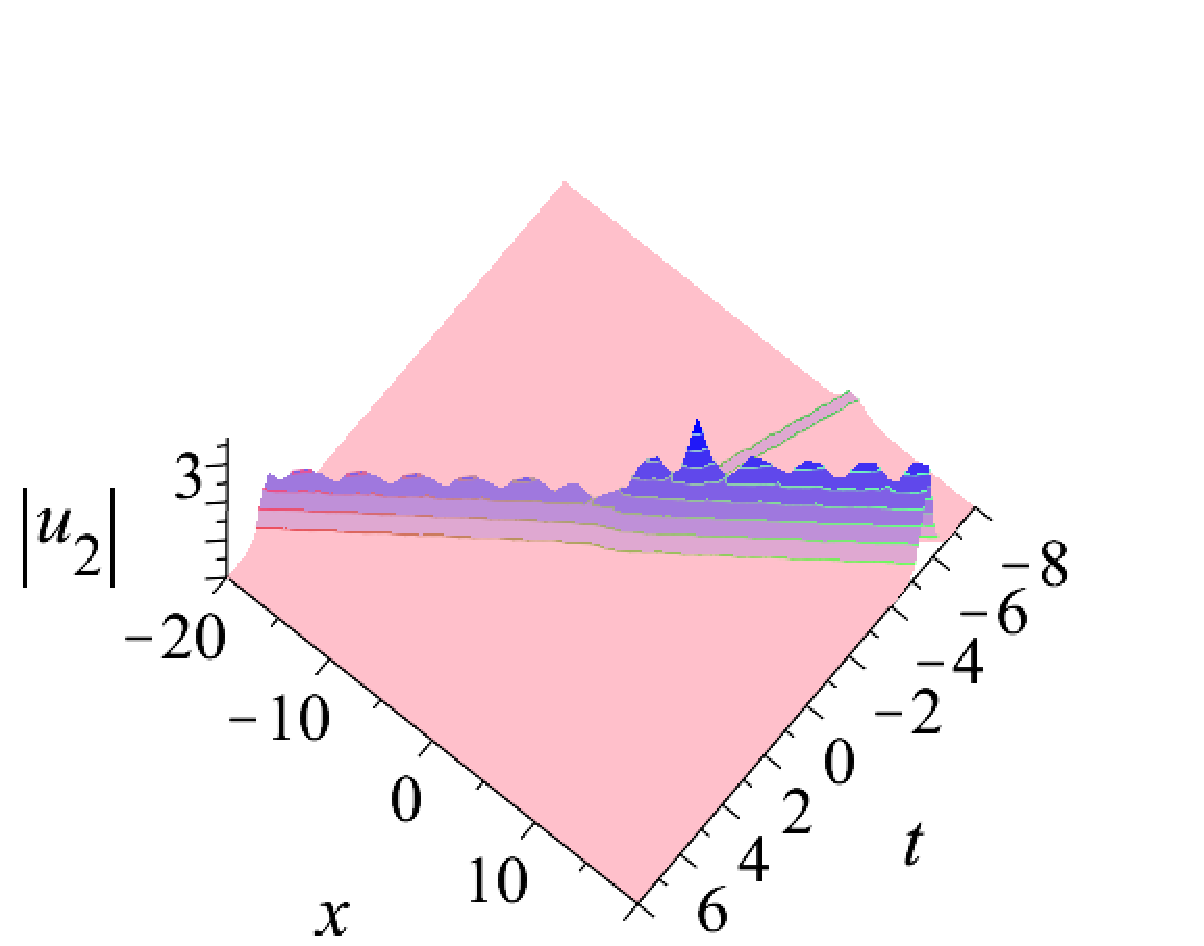}}}
{\rotatebox{0}{\includegraphics[width=4.8cm,height=2.8cm,angle=0]{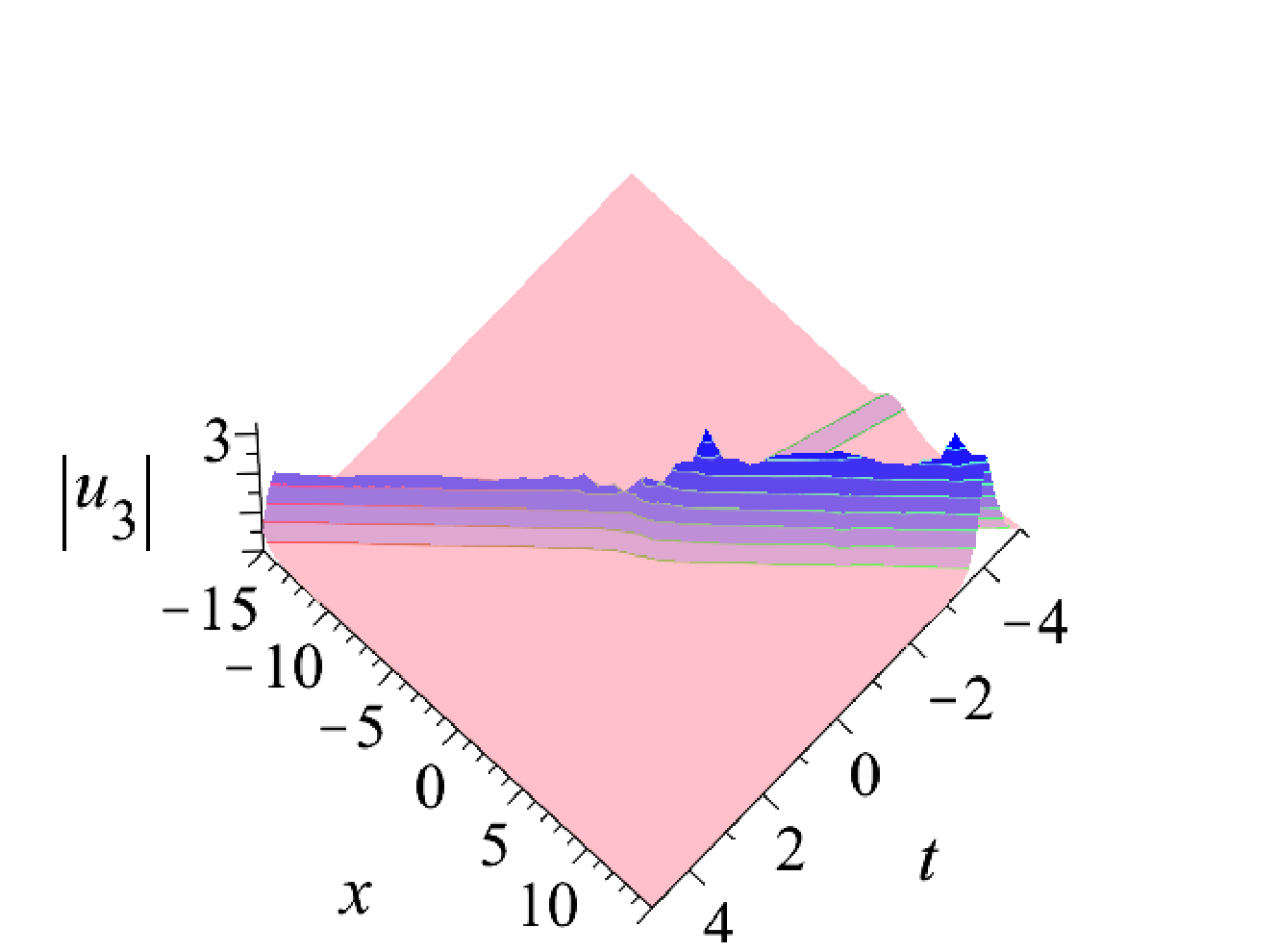}}}

$\qquad\qquad\textbf{(a)}\qquad\qquad\qquad\qquad\qquad\textbf{(b)}
\qquad\qquad\qquad\qquad\qquad\qquad\textbf{(c)}$\\

\noindent
{\rotatebox{0}{\includegraphics[width=4.8cm,height=2.8cm,angle=0]{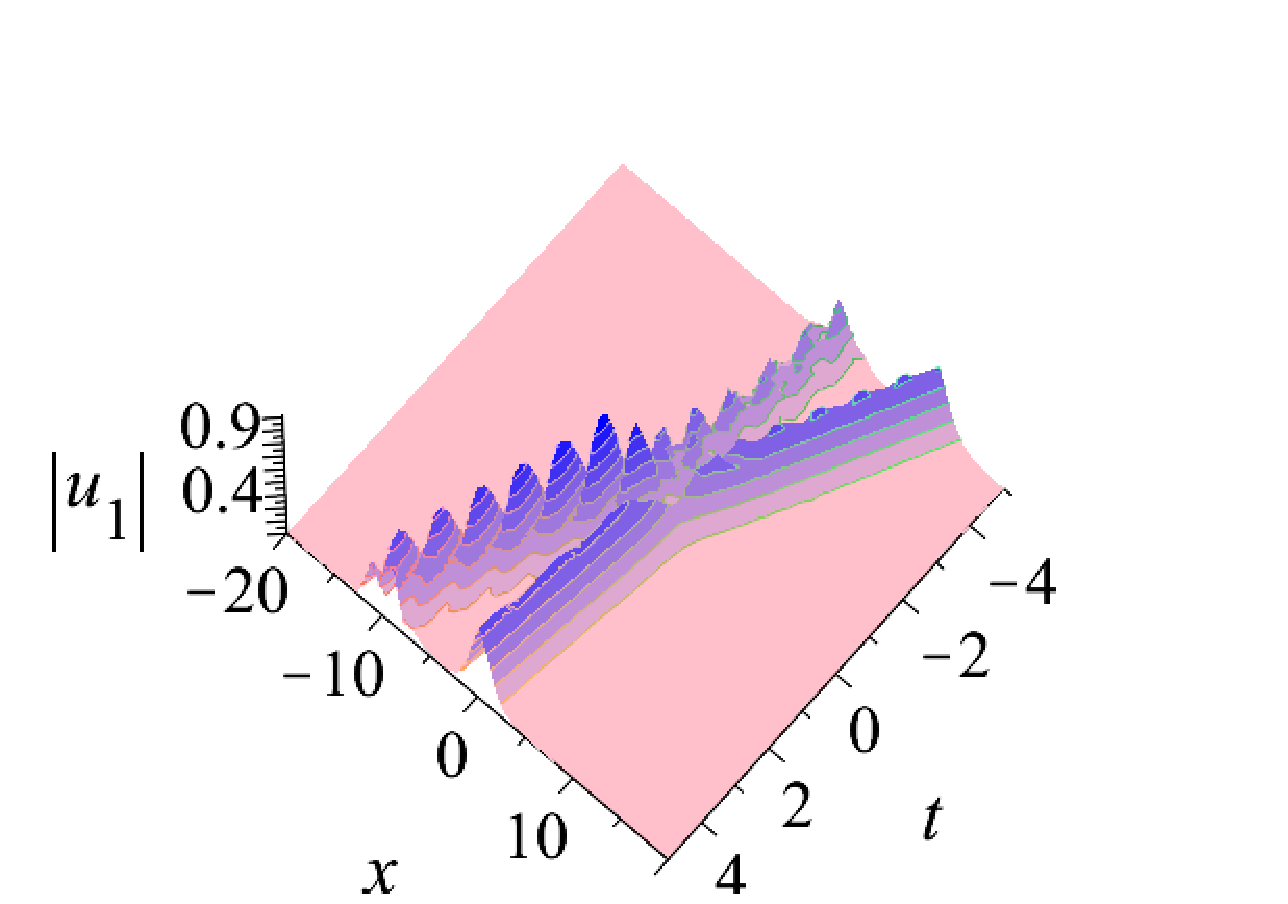}}}
{\rotatebox{0}{\includegraphics[width=4.8cm,height=2.8cm,angle=0]{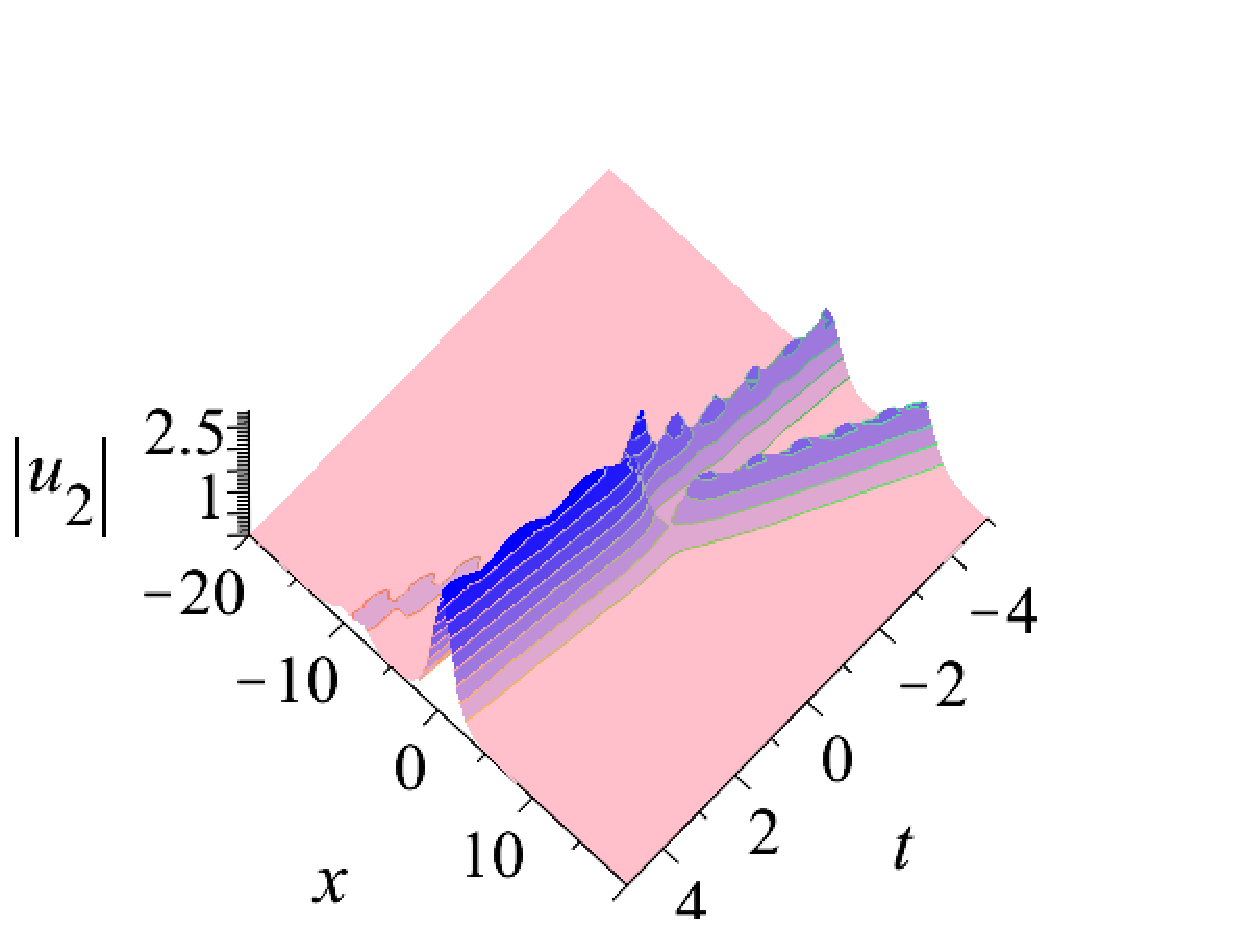}}}
{\rotatebox{0}{\includegraphics[width=4.8cm,height=2.8cm,angle=0]{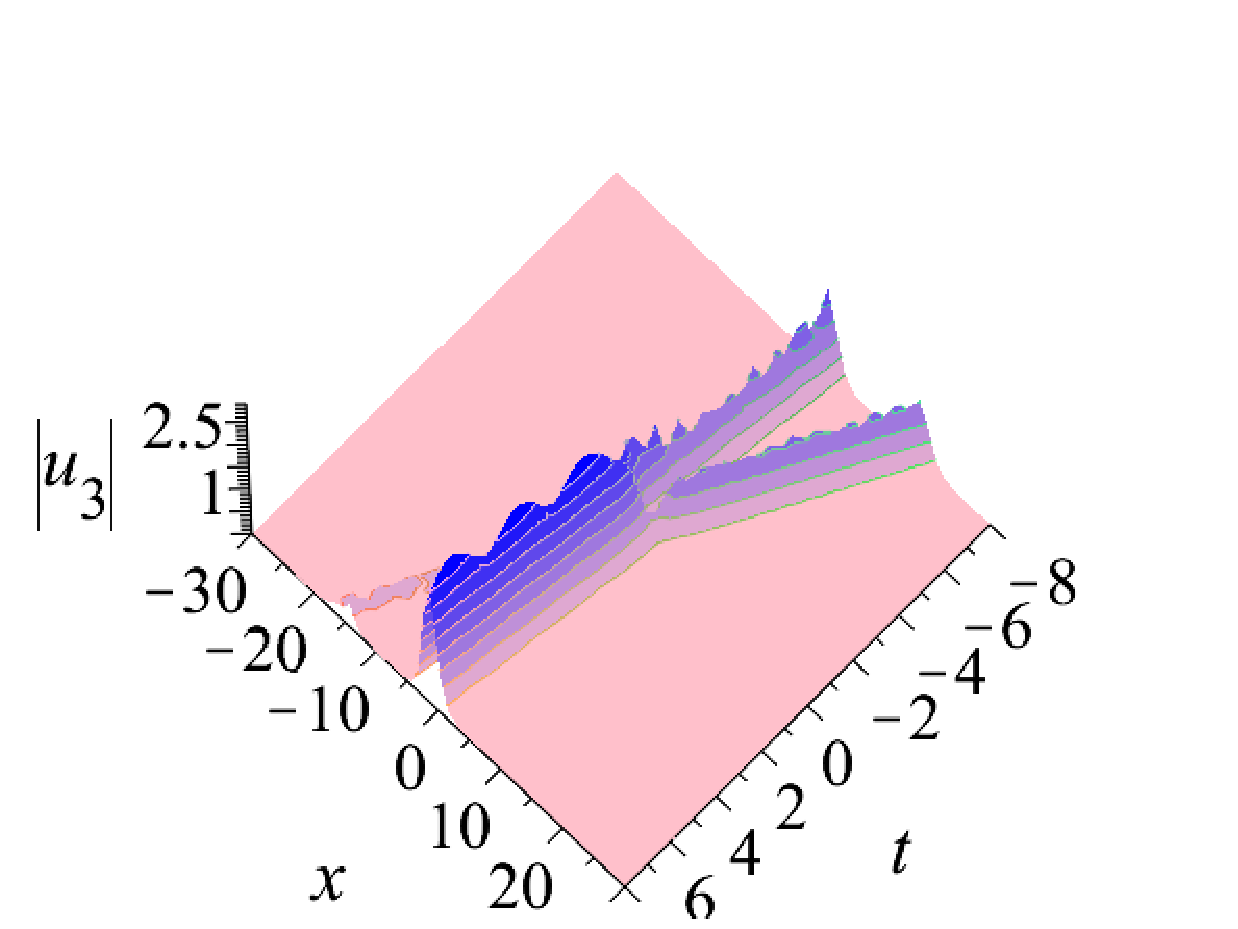}}}

$\qquad\qquad\textbf{(d)}\qquad\qquad\qquad\qquad\qquad\textbf{(e)}
\qquad\qquad\qquad\qquad\qquad\qquad\textbf{(f)}$\\

\noindent { \small \textbf{Figure 2.} (Color online) Two-soliton solutions
with parameters: (\textbf{a},\textbf{b},\textbf{c}):$\lambda_{1}=0.4+0.5i, \lambda_{2}=0.7+0.8i,  \alpha_{1}=\gamma_{1}=\mu_{1}=\rho_{1}=\beta_{2}=\gamma_{2}=\mu_{2}=\rho_{2}=0,\beta_{1}=\alpha_{2}=\gamma_{2}=1$;
(\textbf{d},\textbf{e},\textbf{f}): $\lambda_{1}=0.5+0.5i, \lambda_{2}=0.4+0.6i,  \alpha_{1}=\beta_{1}=\gamma_{1}=\mu_{1}=\rho_{1}=\alpha_{2}=1, \gamma_{2}=2, \beta_{2}=\mu_{2}=\rho_{2}=0$.\\}

\noindent
{\rotatebox{0}{\includegraphics[width=4.8cm,height=2.8cm,angle=0]{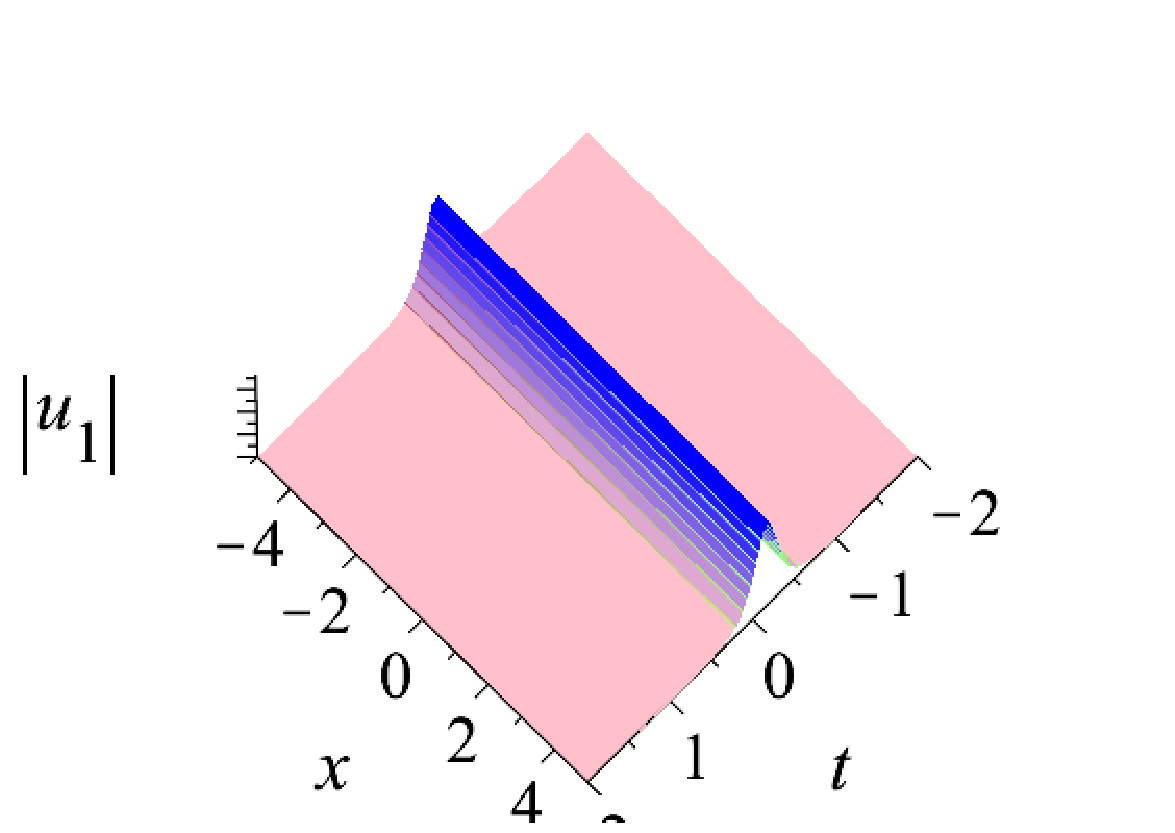}}}
{\rotatebox{0}{\includegraphics[width=4.8cm,height=2.8cm,angle=0]{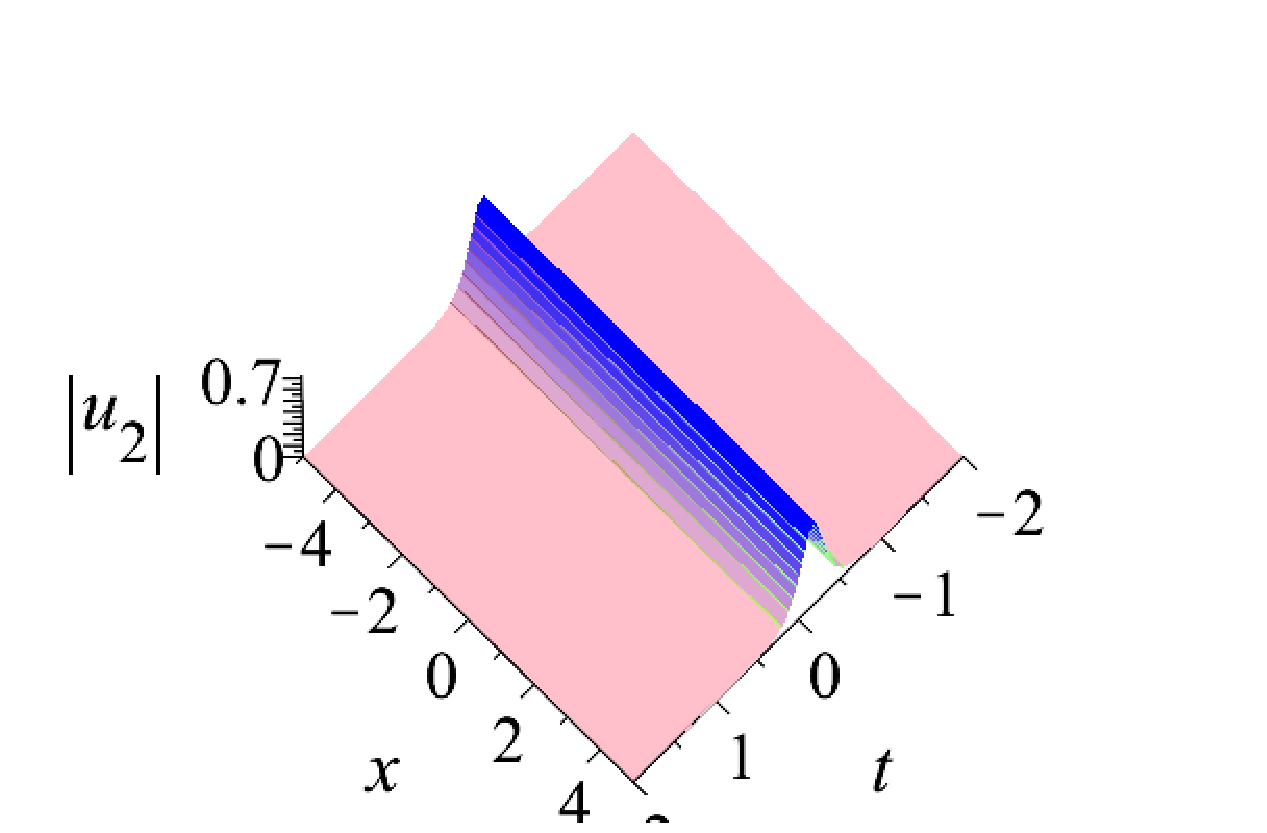}}}
{\rotatebox{0}{\includegraphics[width=4.8cm,height=2.8cm,angle=0]{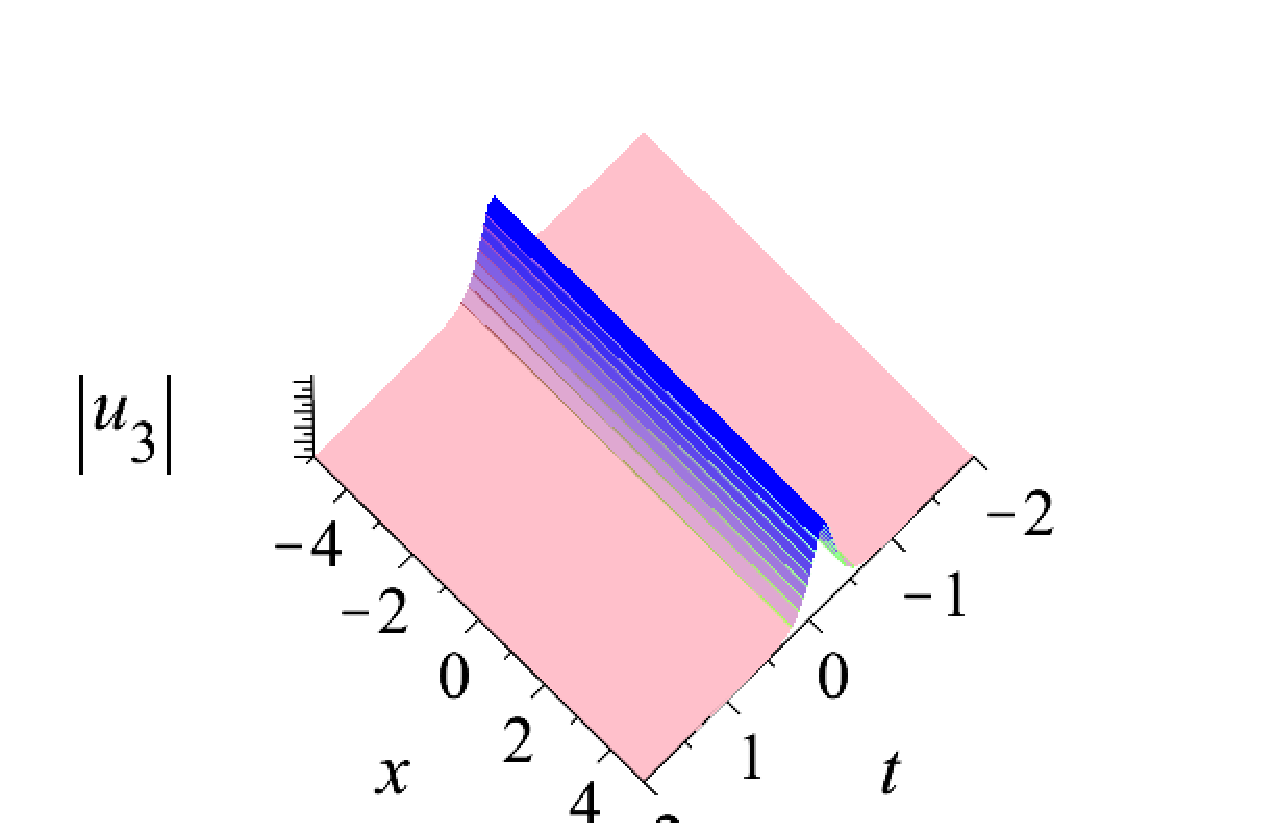}}}

$\qquad\qquad\textbf{(a)}\qquad\qquad\qquad\qquad\qquad\textbf{(b)}
\qquad\qquad\qquad\qquad\qquad\qquad\textbf{(c)}$\\

\noindent
{\rotatebox{0}{\includegraphics[width=4.8cm,height=2.8cm,angle=0]{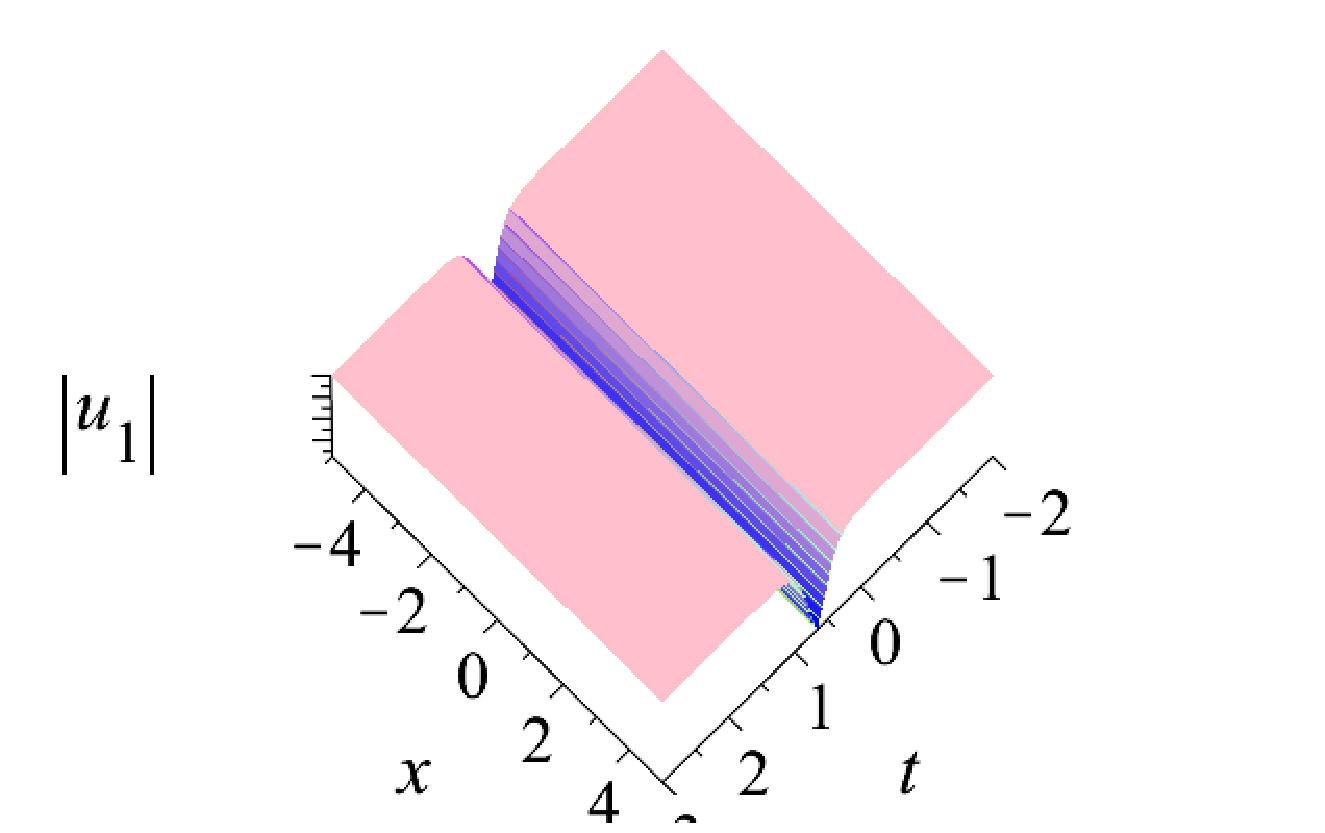}}}
{\rotatebox{0}{\includegraphics[width=4.8cm,height=2.8cm,angle=0]{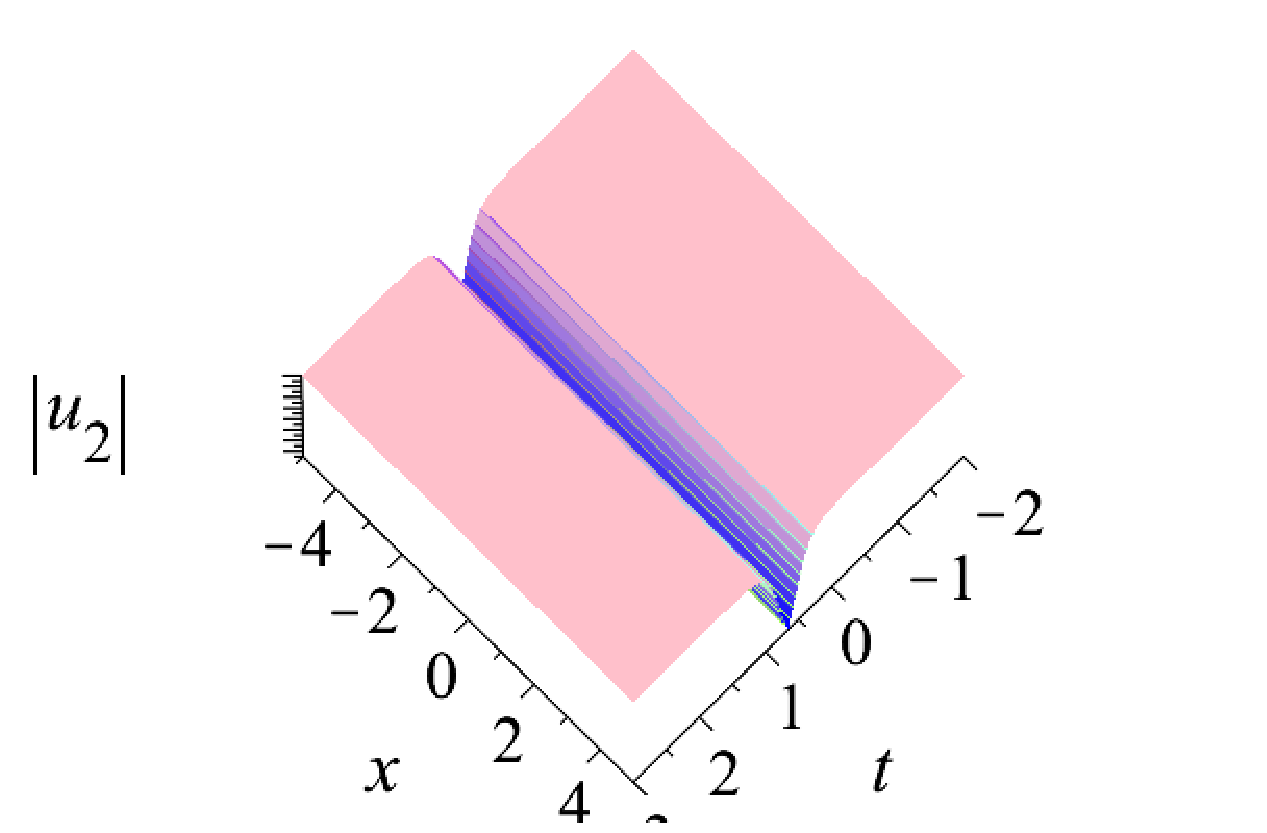}}}
{\rotatebox{0}{\includegraphics[width=4.8cm,height=2.8cm,angle=0]{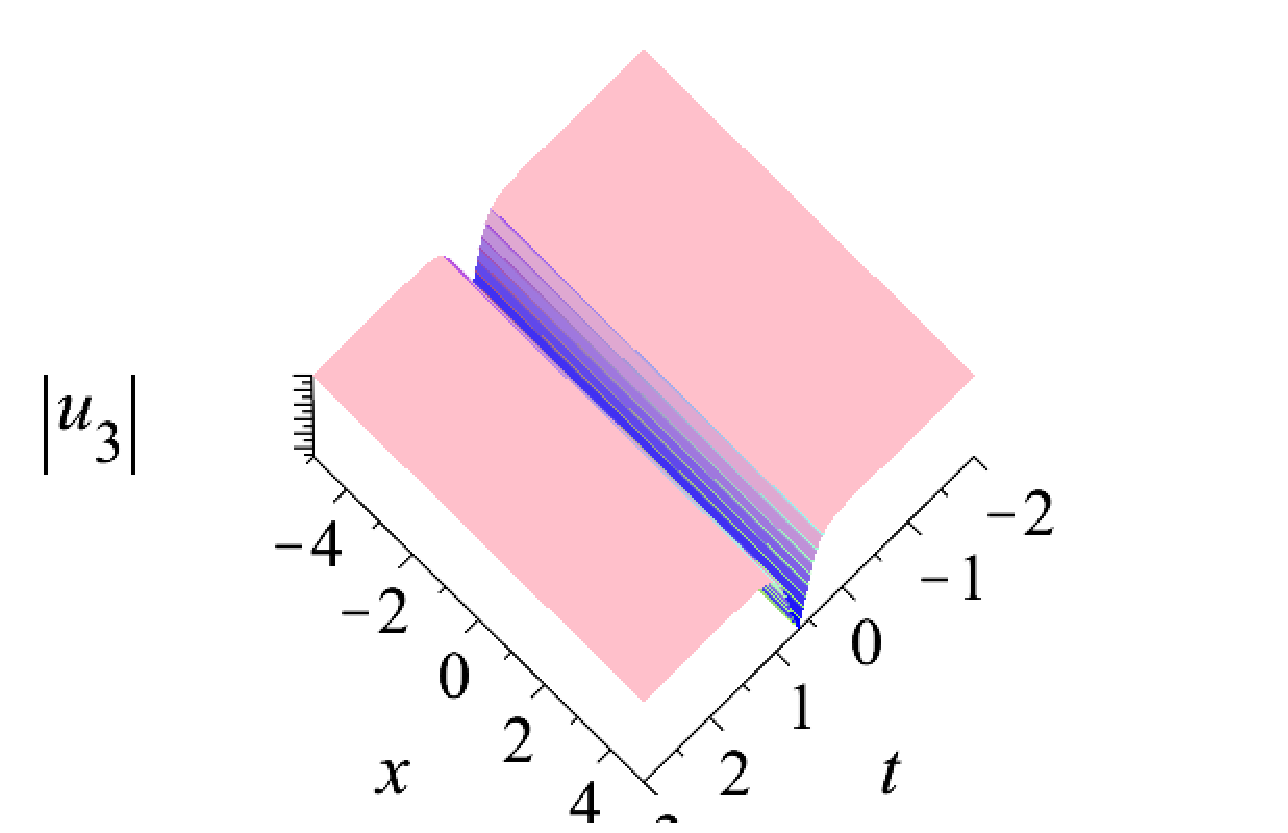}}}

$\qquad\qquad\textbf{(d)}\qquad\qquad\qquad\qquad\qquad\textbf{(e)}
\qquad\qquad\qquad\qquad\qquad\qquad\textbf{(f)}$\\
\noindent { \small \textbf{Figure 3.} (Color online) Bright-dark soliton wave via solutions \eqref{SS-14}
with parameters: $\alpha_{1}=1, \gamma_{1}=2, \rho_{1}=3$.
(\textbf{a},\textbf{b},\textbf{c}): $\eta_{1}=1$. (\textbf{d},\textbf{e},\textbf{f}): $\eta_{1}=-1$.\\}

\noindent
{\rotatebox{0}{\includegraphics[width=4.8cm,height=2.8cm,angle=0]{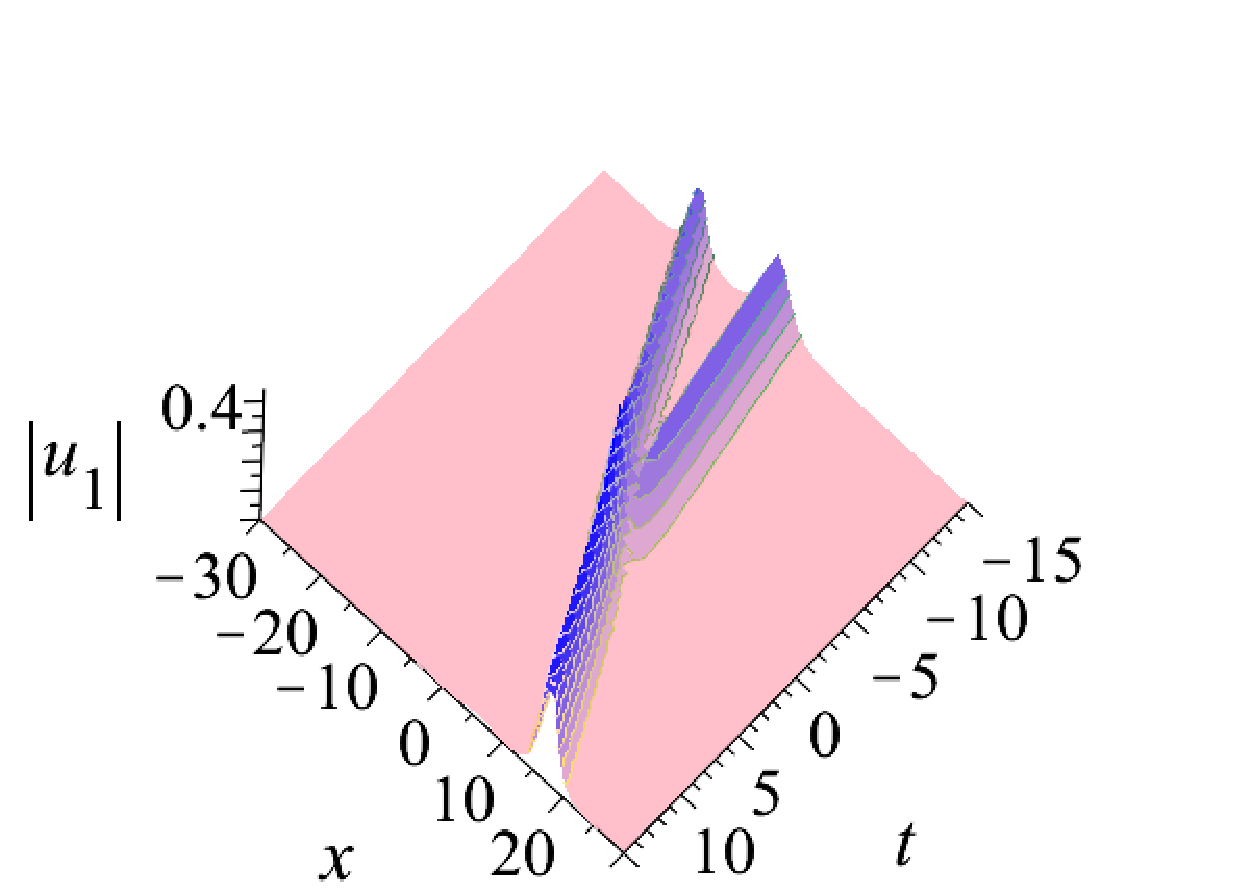}}}
{\rotatebox{0}{\includegraphics[width=4.8cm,height=2.8cm,angle=0]{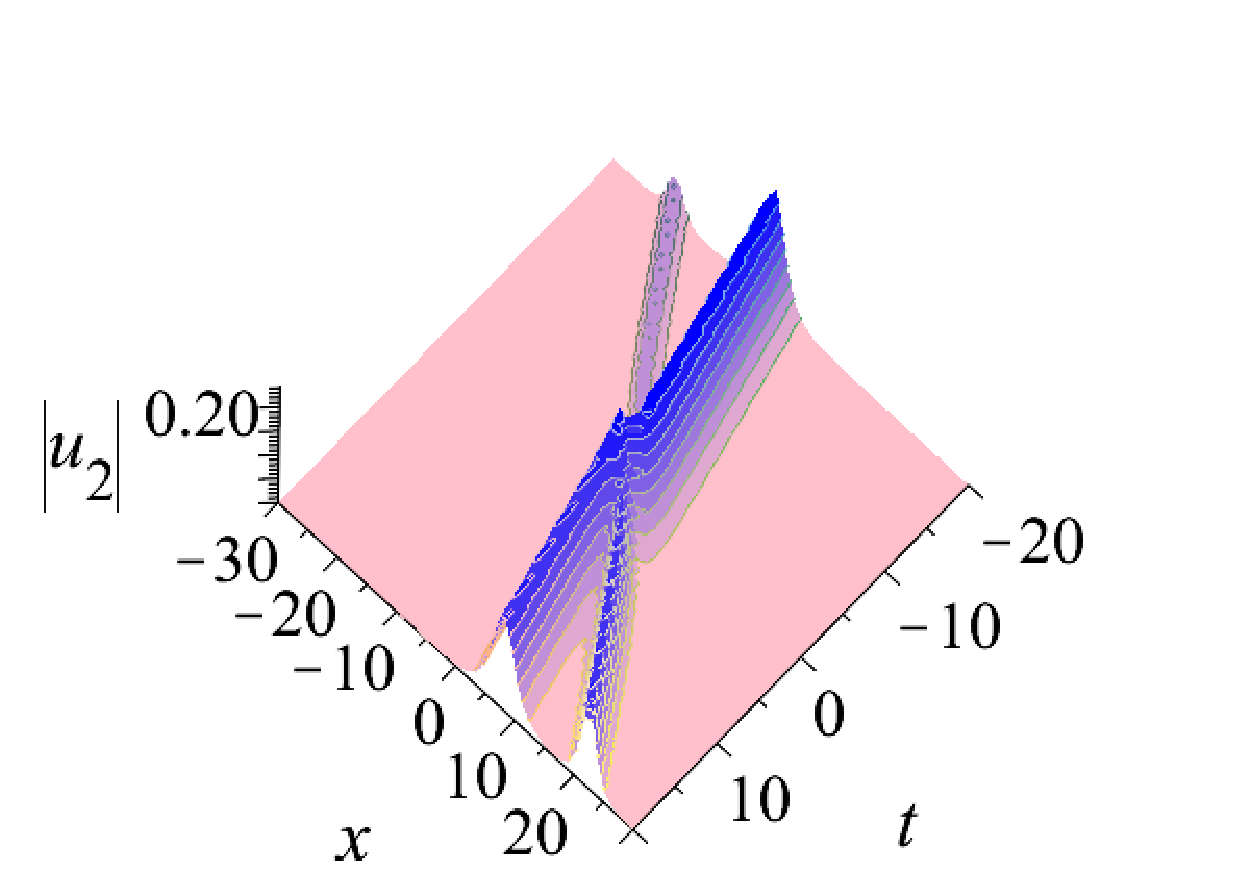}}}
{\rotatebox{0}{\includegraphics[width=4.8cm,height=2.8cm,angle=0]{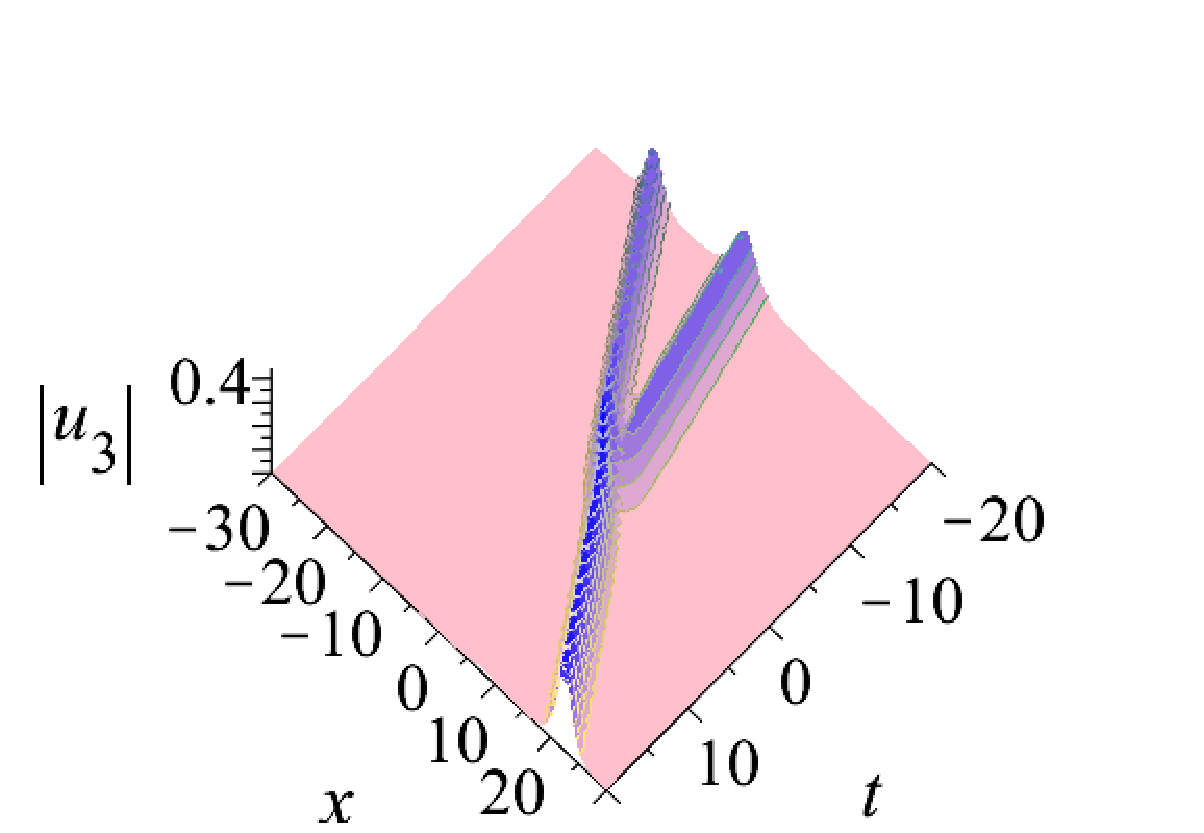}}}

$\qquad\qquad\textbf{(a)}\qquad\qquad\qquad\qquad\qquad\textbf{(b)}
\qquad\qquad\qquad\qquad\qquad\qquad\textbf{(c)}$\\

\noindent
{\rotatebox{0}{\includegraphics[width=4.8cm,height=2.8cm,angle=0]{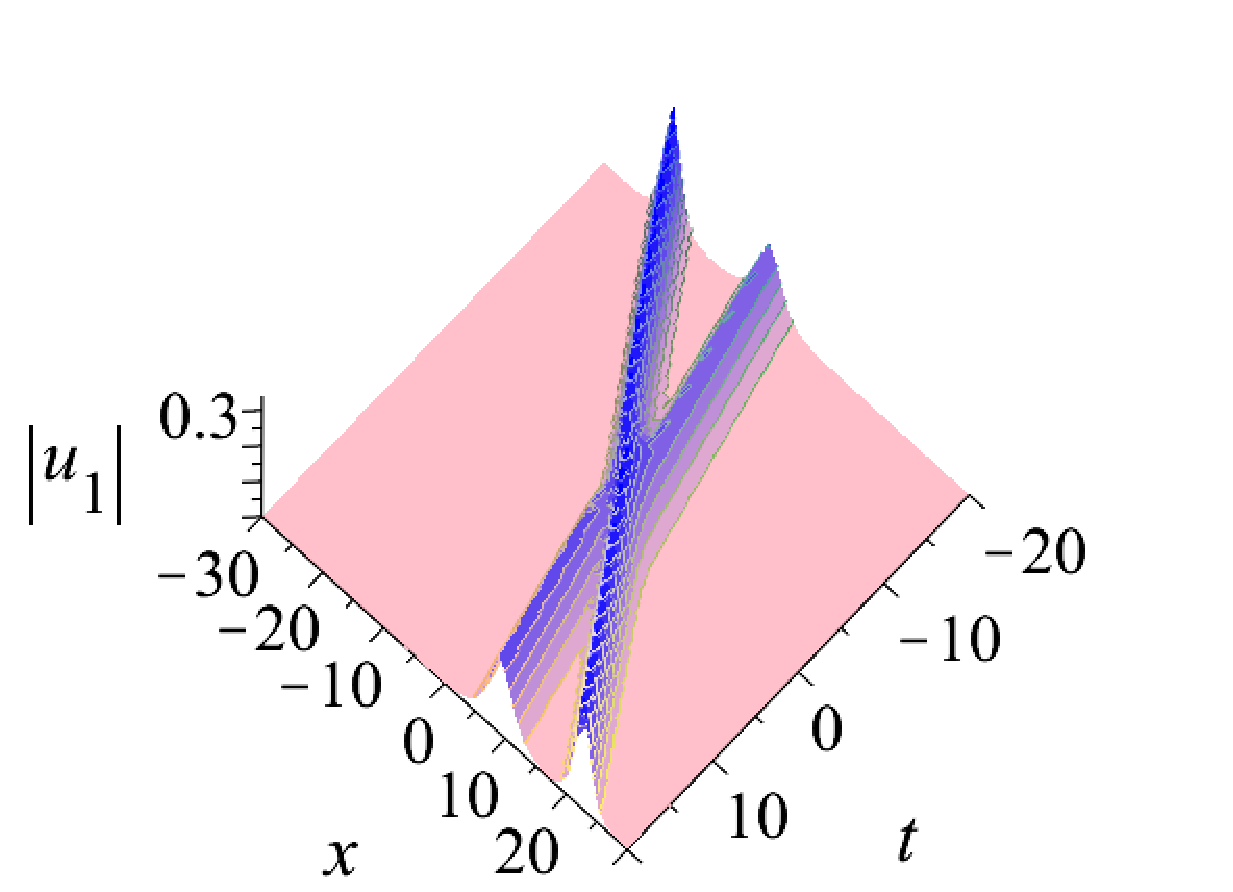}}}
{\rotatebox{0}{\includegraphics[width=4.8cm,height=2.8cm,angle=0]{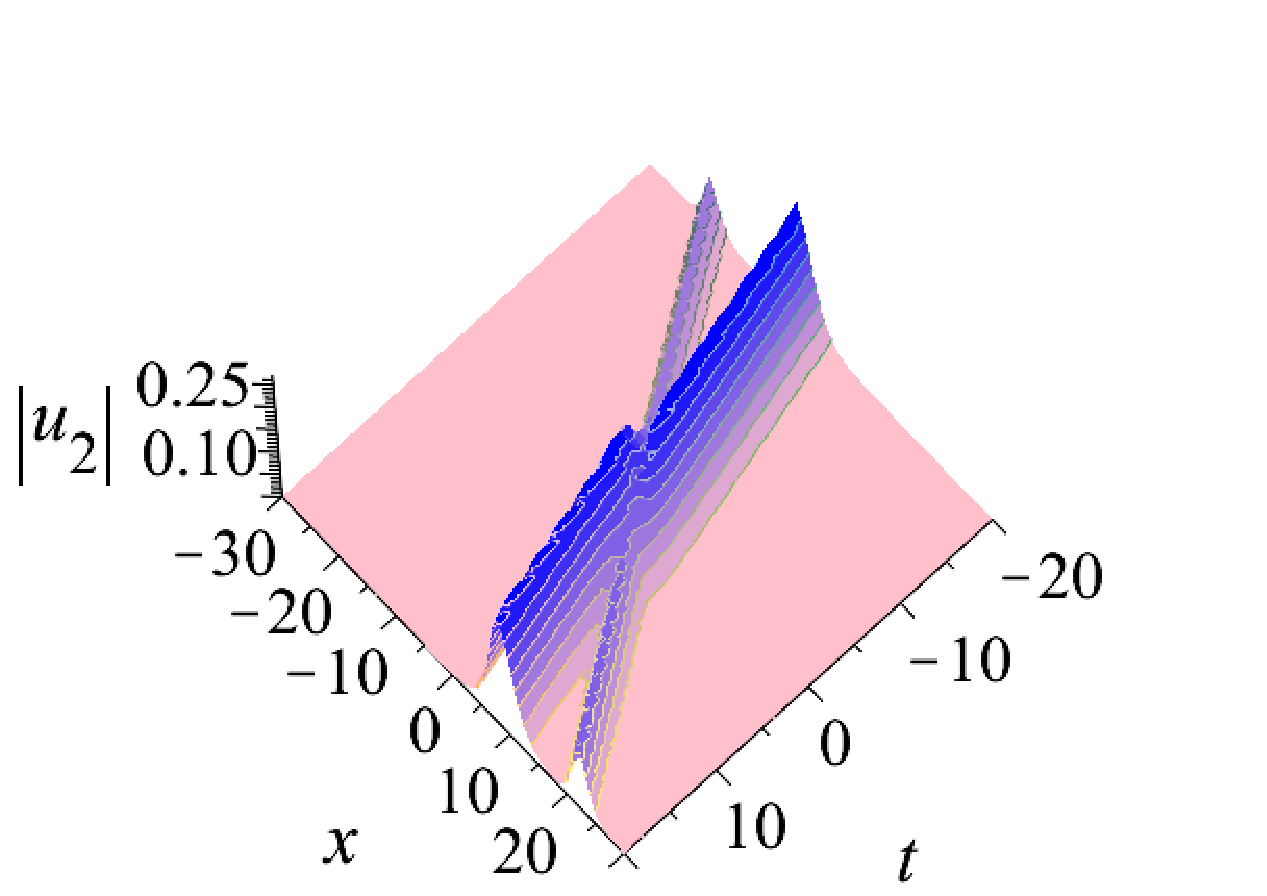}}}
{\rotatebox{0}{\includegraphics[width=4.8cm,height=2.8cm,angle=0]{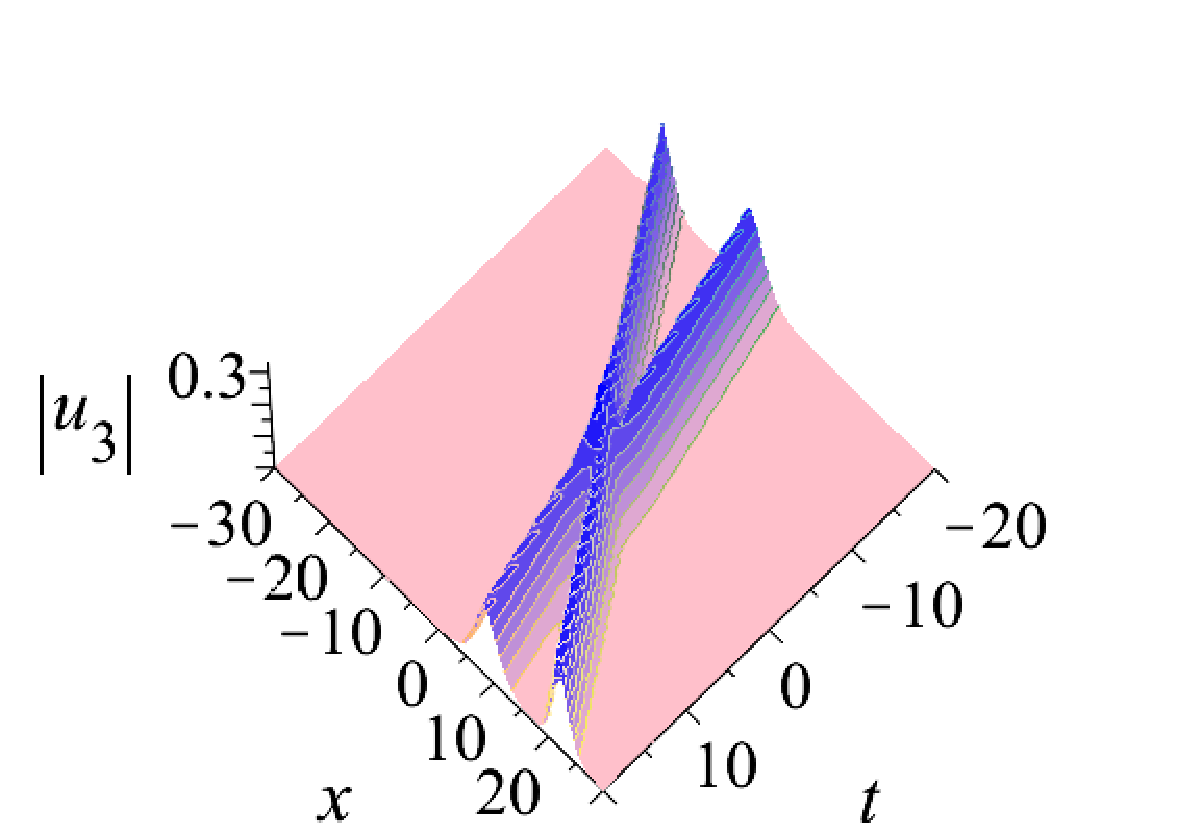}}}

$\qquad\qquad\textbf{(d)}\qquad\qquad\qquad\qquad\qquad\textbf{(e)}
\qquad\qquad\qquad\qquad\qquad\qquad\textbf{(f)}$\\

\noindent { \small \textbf{Figure 4.} (Color online) Two-soliton solutions \eqref{SS-12} (N=2)
with parameters: (\textbf{a},\textbf{b},\textbf{c}): $\lambda_{1}=0.3i, \lambda_{2}=0.5i,  \alpha_{1}=\alpha_{2}=\gamma_{1}=\rho_{1}=\rho_{2}=1, \gamma_{2}=0.5$;
(\textbf{d},\textbf{e},\textbf{f}):$\lambda_{1}=0.3i, \lambda_{2}=0.5i,  \alpha_{1}=1, \alpha_{2}=i, \gamma_{1}=1+i, \rho_{1}=1+i, \rho_{2}=i, \gamma_{2}=0.5i$.\\}

%\noindent
%{\rotatebox{0}{\includegraphics[width=4.5cm,height=3.9cm,angle=0]{6-1.eps}}}
%~~~
%{\rotatebox{0}{\includegraphics[width=4.5cm,height=3.9cm,angle=0]{6-2.eps}}}
%~~~
%{\rotatebox{0}{\includegraphics[width=4.5cm,height=3.9cm,angle=0]{6-3.eps}}}
%
%$\qquad\qquad\textbf{(a)}\qquad\qquad\qquad\qquad\qquad\qquad\textbf{(b)}
%\qquad\qquad\qquad\qquad\qquad\qquad\textbf{(c)}$\\
%
%\noindent { \small \textbf{Figure 6.} (Color online) Two-soliton solutions \eqref{SS-12} (N=2)
%with parameters: $\lambda_{1}=0.3i, \lambda_{2}=0.5i,  \alpha_{1}=1, \alpha_{2}=i, \gamma_{1}=1+i, \rho_{1}=1+i, \rho_{2}=i, \gamma_{2}=0.5i$.\\}

\section{Conclusions}

This work is dedicated to the associated multi-soliton solutions of the three-component coupled Sasa-Satsuma equation \eqref{SSE}.
By using the Riemann-Hilbert method, we have obtained a class of multi-soliton solutions for the Eq. \eqref{SSE}.
Additionally, based on the soliton solution formulas, we find some interesting soliton solutions which contain
breather-type solution, single-soliton solution etc.
In order to help readers understand those soliton solutions better,
the propagation behaviors have been shown by graphical simulations (i.e.,Figs.1-4).
Particularly, the Eq. \eqref{SSE}
we investigated in this work are fairly more general as they involve a $7\times7$ Lax pair.
The celebrated Sasa-Satsuma equation, a crucial model in fiber optics, is its particular case.
Another important reduction of the Eq. is the coupled Sasa-Satsuma equation.
Consequently, the $N$-soliton solutions of the celebrated Sasa-Satsuma equation and
the coupled Sasa-Satsuma equation \cite{gxg-cnsns} can be respectively obtained by reducing the $N$-soliton solutions of \eqref{SSE}.
More importantly,
in certain physical situations, two or more wave packets of different carrier frequencies
appear simultaneously, and their interactions can be governed by the coupled NLS equations.
Examples include nonlinear light propagation in a birefringent optical fiber or a
wavelength-division-multiplexed system \cite{zjj-1,zjj-2,zjj-3}, spinor Bose-Einstein condensates (BECs) \cite{zjj-6,zjj-7}, the
interaction of Bloch-wave packets in a periodic system \cite{zjj-4}, the evolution of two surface wave packets in deep water \cite{zjj-5}, etc.
Since the coupled
NLS equations arise in a wide variety of physical subjects
such as nonlinear optics, water waves, BECs, etc,
these results should prove useful to the investigations of those physical problems.

Finally, we remark that there are several methods to get exact solutions for nonlinear evolution equations (NLEEs),
such as the
Darboux transformation \cite{vbm-1991,mwx-LMP}, the RH problem approach \cite{SP-1984,MJ-1991},
the Hirota method \cite{RH-2004}, the dressing method \cite{dm-1}, the Wronskian technique \cite{wr-1}, etc.
Thus, it is very necessary to discuss whether the Eq. \eqref{SSE} can be solved by using these approaches?
These will be left for future discussions.
Recently, the Darboux transform method is used to derive the $N$-soliton solutions of the celebrated Sasa-Satsuma equation.
Comparing with the soliton solution formulae obtained in this work and those constructed by Darboux
transformation in Ref. \cite{llm-1}, it is very clear that Eq.\eqref{SS-12} is much simpler.
In Ref. \cite{hrt-1}, the bilinear method is used to derive the soliton solutions of the integrable Sasa-Satsuma equation.
Unfortunately, the construction of multisoliton solutions to this equation
presents difficulties due to its complicated bilinearization.
So the authors discuss briefly some previous attempts and then
present the correct bilinearization based on the interpretation of the SS equation as a reduction of the
three-component KP hierarchy.
As a result, we find that the RH method can provide an effective
and powerful mathematical tool to derive exact solutions of NLEEs, which should be suitable to
analyze other models in mathematical physics and engineering.

\section*{Acknowledgements}
\hspace{0.2cm}%We express our sincere thanks to all the persons who have provided valuable suggestions to this paper.
This work is supported
by the National Key Research and Development Program of China under Grant No. 2017YFB0202901 and the National Natural Science Foundation of China under Grant No. 11871180.

%\end{CJK*}
%\end{document}

% Converted from Microsoft Word to LaTeX
% by Chikrii SoftLab Word2TeX converter (version 2.4)
% Copyright (C) 1999-2001 Kirill A. Chikrii, Anna V. Chikrii
% Copyright (C) 1999-2001 Chikrii SoftLab.
% All rights reserved.
% http://www.word2tex.com/
% mailto: info@word2tex.com, support@word2tex.com

% Warning: You are using UNREGISTERED Chikrii SoftLab Word2TeX!
%          In UNREGISTERED mode some restrictions will apply.
%          For more information please visit http://www.word2tex.com/
% YOU CAN USE THIS FILE WITH THE SOLE PURPOSE OF EVALUATING Word2TeX.

%\documentclass [12pt]{article}

%\begin{document}

\end{document}